\documentclass[12pt,reqno]{amsart}

\usepackage{latexsym}
\usepackage{amssymb}

\renewcommand{\Re}{{\operatorname{Re}\,}}
\renewcommand{\Im}{{\operatorname{Im}\,}}

\renewcommand{\epsilon}{\varepsilon}
\newcommand{\dist}{{\operatorname{dist}}}
\newcommand{\var}{{\operatorname{Var}}}

\newcommand{\sm}{\smallsetminus}
\newcommand{\szego}{Szeg\"o }

\newcommand{\inv}{^{-1}}
\newcommand{\kahler}{K\"ahler }
\newcommand{\sqrtn}{\sqrt{N}}
\newcommand{\wt}{\widetilde}
\newcommand{\wh}{\widehat}
\newcommand{\PP}{{\mathbb P}}
\newcommand{\N}{{\mathbb N}}
\newcommand{\R}{{\mathbb R}}
\newcommand{\C}{{\mathbb C}}

\newcommand{\Z}{{\mathbb Z}}

\newcommand{\CP}{\C\PP}
\renewcommand{\d}{\partial}
\newcommand{\dbar}{\bar\partial}
\newcommand{\ddbar}{\partial\dbar}
\newcommand{\U}{{\rm U}}

\newcommand{\E}{{\mathbf E}}
\renewcommand{\H}{{\mathbf H}}

\newcommand{\half}{{\frac{1}{2}}}
\newcommand{\vol}{{\operatorname{Vol}}}

\newcommand{\SU}{{\operatorname{SU}}}
\newcommand{\FS}{{{\operatorname{FS}}}}

\renewcommand{\phi}{\varphi}

\newcommand{\ccal}{\mathcal{C}}
\newcommand{\dcal}{\mathcal{D}}

\newcommand{\fcal}{\mathcal{F}}

\newcommand{\hcal}{\mathcal{H}}
\newcommand{\ical}{\mathcal{I}}

\newcommand{\lcal}{\mathcal{L}}

\newcommand{\ncal}{\mathcal{N}}
\newcommand{\ocal}{\mathcal{O}}
\newcommand{\pcal}{\mathcal{P}}

\newcommand{\scal}{\mathcal{S}}

\newcommand{\al}{\alpha}

\newcommand{\ga}{\gamma}

\newcommand{\La}{\Lambda}
\newcommand{\la}{\lambda}
\newcommand{\ep}{\varepsilon}
\newcommand{\de}{\delta}

\newcommand{\om}{\omega}
\newcommand{\Om}{\Omega}

\newtheorem{theo}{{\sc Theorem}}[section]

\newtheorem{cor}[theo]{{\sc Corollary}}

\newtheorem{lem}[theo]{{\sc Lemma}}
\newtheorem{prop}[theo]{{\sc Proposition}}

\newenvironment{rem}{\medskip\noindent{\it Remark:\/} }{\medskip}
\newtheorem{defin}[theo]{{\sc Definition}}

\title
{Number variance of random zeros}

\author{Bernard Shiffman}
\author{Steve Zelditch}
\address{Department of Mathematics, Johns Hopkins University, Baltimore, MD
21218, USA} \email{bshiffman@jhu.edu, szelditch@jhu.edu}

\thanks{Research partially supported by NSF grant
DMS-0100474 (first author) and  DMS-0302518 (second author).}

\date{January 6, 2006}

\begin{document}

\begin{abstract}

The main results of this article are  asymptotic formulas for the
variance of the number of zeros of a Gaussian  random
polynomial of degree $N$ in an open set $U \subset \C$ as the
degree $N \to \infty$, and more generally for the zeros of random
holomorphic sections  of high powers of any
positive line bundle over any Riemann surface. The formulas were
conjectured in special cases by  Forrester and Honner.  In
higher dimensions, we give similar formulas for the variance of
the volume inside a domain $U$ of the zero hypersurface
of a random holomorphic section of a high power of a
positive line bundle over any compact K\"ahler manifold. These results
generalize the variance asymptotics of Sodin and Tsirelson
for special model ensembles of  chaotic analytic functions in one
variable to any ample line bundle and Riemann surface. We also combine
our methods with those of Sodin-Tsirelson to generalize their
asymptotic normality results for smoothed number statistics.
\end{abstract}

\maketitle

\tableofcontents

 \section{Introduction}

This article is concerned with number and volume statistics for
Gaussian random holomorphic functions (and sections). To introduce
our subject, let us start with the simplest case of  holomorphic
polynomials $p_N$ of degree $N$  of one complex variable.  By
homogenizing, we may identify the space $\pcal_N$ of polynomials
of degree $N$ with the space $H^0(\CP^1, \ocal(N))$ of holomorphic
sections of the $N$-th power of the hyperplane section bundle over
$\CP^1$. This space carries a natural $SU(2)$-invariant inner
product and associated Gaussian measure $\gamma_N$. To each
polynomial $p_N$ we associate its zero set $Z_{p_N} \subset \CP^1$
and thus obtain a random point process on $\CP^1$. Given an open
subset $U \subset \CP^1$, we define the integer-valued  random
variable
\begin{equation} \ncal_N^U(p_N) = \#\{z\in U:  p_N(z)=0\}
\end{equation}
on $\pcal_N$ counting the number of zeros of $p_N$ which lie in
$U$. Clearly, $\ncal^U_N$  is discontinuous along the set of polynomials
having a zero on the boundary $\partial U$. It is easy to see from the
$\SU(2)$ invariance that the expected value of this random
variable is given by
\begin{equation*}\E (\ncal^U_N) = N\int_U \frac
i{2\pi}\,\Theta_h\;,\end{equation*} where
$\E(X)$ denotes the expectation of a random variable $X$ and where
$\Theta_h$ is the curvature form of the Fubini-Study metric; i.e.,
the expected zero distribution is uniform on $\CP^1$ with respect to
its
$SU(2)$ invariant area form. The variance
$$\var(\ncal^U_N) = \E \big( \ncal^U_N - \E(\ncal^U_N) \big)^2$$
of $\ncal^U_N$ measures  the fluctuations of $\ncal^U_N$, i.e.\ 
the extent to
which the number of zeros of individual polynomials conforms to or
deviates from  the expected number. More generally, we can study
the same problem for Gaussian random holomorphic sections $s_N \in
H^0(M, L^N)$ of powers of any positive holomorphic line bundle $L
\to M$ over any compact Riemann surface $M$. In  \cite{SZ}, we
showed that in this case, the expected value of the random variable
$\ncal^U_N$ has the asymptotics
\begin{equation}\label{MEANa}\frac 1
{N}\,\E (\ncal^U_N) = \frac i{2\pi}\int_U \,\Theta_h\ +\ O\left(\frac
1N\right)\;.\end{equation}

Our first result  gives an  estimate for the variance of the
number of zeros on a domain in a compact complex curve, extending
and sharpening a result of Forrester and Honner \cite{FH} (see
also Sodin-Tsirelson \cite{ST}):

\begin{theo}\label{number} Let $(L,h)$ be a
positive Hermitian holomorphic line bundle over a compact complex
curve $M$. We give $H^0(M,L^N)$ the Hermitian Gaussian measure
induced by $h$ and the area form $\om= \frac i2 \Theta_h$. Let $U$
be a domain in $M$ with piecewise $\ccal^2$ boundary and no cusps.
Then for random sections $s_N\in H^0(M,L^N)$, we have
$$\var\big(\#\{z\in U:s_N(z)=0\}\big) =
\sqrtn\left[\frac{\zeta(3/2)}{8\pi^{3/2}}\, \mbox{\rm Length}(\d
U)
 +O(N^{-\half +\ep})\right]\;.$$
\end{theo}

This theorem  proves  a strong form of self-averaging for the
number of zeros in $U$. Here,  a sequence $X_N$ of random
variables is called {\it self-averaging\/} if the fluctuations of
$X_N$ are of smaller order than its typical values, or in other
words if $\frac{ \var(X)}{ (\E X)^2} \to 0$.

In higher dimensions, the analogous point process is
defined by the simultaneous zeros of $m$ polynomials 
in $m$ variables, or more generally, $m$ sections on an $m$-dimensional complex
manifold. In this article we  consider instead the simpler `volume analogue' of number
statistics for one random polynomial or section $s$ in  $m$ dimensions. We let $Z_{s}$
denote the zero set of the random holomorphic section $s$. Recall that the volume of
$Z_{s}$ in a domain $U$ is given by
$$\vol_{2m-2}[Z_{s}\cap U]=\int_{Z_{s}\cap U} \frac
1{(m-1)!}\om^{m-1}=
\left([Z_{s}],\ \chi_U\,\frac 1{(m-1)!}\om^{m-1}\right)\;,$$ where
$[Z_{s}]$ denotes the current of integration over $Z_s$.  Our higher
dimensional generalization of Theorem
\ref{number} is the following asymptotic   formula for the variance of the
volume of the zero divisor in a domain with nice boundary.

\begin{theo}\label{volume}  Let $(L,h)$ be a
positive Hermitian holomorphic line bundle over a compact \kahler
manifold $(M,\om)$, where $\om= \frac i2 \Theta_h$. We give $H^0(M,L^N)$ the
Hermitian Gaussian measure induced by $h,\om$ (see
Definition~\ref{HG}).  Let $U$ be a domain in $M$ with piecewise
$\ccal^2$ boundary and no cusps.  Then for  random sections
$s_N\in H^0(M,L^N)$, we have
$$\var\big(\vol_{2m-2}[Z_{s_N}\cap U]\big) =
N^{-m+3/2}\left[\nu_m\,\vol_{2m-1}(\d U)
 +O(N^{-\half +\ep})\right]\;,$$ where $$\nu_m =
\frac{\pi^{m-5/2}}{8}\;\zeta(m+\textstyle \half)\;.$$
\end{theo}

Here, we say that $U$ has piecewise $\ccal^k$ boundary without
cusps if for each boundary point $z_0\in\d U$, there exists a
(not necessarily convex) closed polyhedral cone
$C\subset\R^{2m}$ and a $\ccal^k$ diffeomorphism
$\rho:V\to
\rho(V)\subset\R^{2m}$, where
$V$ is a neighborhood of $z_0$, such that $\rho (V\cap
\overline U) = \rho(V)\cap C$.  If $M$ is a complex curve, this
condition  means that $\d U$ is a piecewise $\ccal^k$ curve
with distinct tangents at corners and self-intersection
points.

A model case of Theorem \ref{volume} (as well as of the results stated
below) is where $M=\CP^m$ and $L=\ocal(1)$ is the hyperplane section bundle with the
$\SU(m+1)$-invariant Hermitian metric $h$.  Then sections in $H^0(M,L^N)$ are
homogeneous degree $N$ holomorphic polynomials on $\C^{m+1}$, and volumes are computed
with respect to the Fubini-Study metric $\om_\FS= \frac i2 \Theta_h = \frac i2 \ddbar
\log |z|^2$ on $\CP^m$.

In addition to the number variance problem raised by
Forrester-Honner \cite{FH}, the main motivation for this article came
from the variance and asymptotic normality theorems of Sodin-Tsirelson
\cite{ST}  for certain model  random analytic functions. They
consider the  smooth analogue of number statistics, sometimes
called  `linear statistics', defined by the random variables
\begin{equation} N_{\phi} (s) =(Z_s,\phi)= \sum_{\{z\in M:s(z)=0\}}
\phi(z),  \end{equation} where $\phi \in \ccal^3_c(M)$ is a test
function.  In our early paper \cite{SZ}, we  showed that
\begin{equation}\label{ave} \E(Z_{s_N},\phi)= N \int_M \om\wedge\phi
+O(1)\;,\end{equation} and we gave a crude
bound (see \cite[Lemma 3.3]{SZ})
\begin{equation}\label{polydiv}\frac{\var (Z_N,\phi)}{\big[\, \E(Z_N,\phi) \,\big]^2} = O\left(\frac
1{N^2} \right)\end{equation} on the variance, which was sufficient
to prove a strong law of large numbers for the distribution of
zeros. In certain model ensembles,  Sodin-Tsirelson \cite{ST}
improved this result to a sharp estimate as an ingredient in their
asymptotic normality result for zeros. Our next result generalizes
their variance asymptotics  for the zeros of random polynomials
$f_N$ of degree $N$ (and their counterparts for model chaotic
analytic functions in $\ocal(D)$ and $\ocal(\C)$) to any compact
\kahler manifold:

\begin{theo}\label{sharp} Let $(L,h)$ be a
positive Hermitian holomorphic line bundle over a compact \kahler
manifold $(M,\om)$, where $\om= \frac i2 \Theta_h$. Let
$\phi$ be a real $(2m-2)$-form on $M$ with $\ccal^3$ coefficients.  We
give
$H^0(M,L^N)$ the Hermitian Gaussian measure induced by $h$ and the area form
$\om$ (see Definition~\ref{HG}).  Then for  random sections $s_N\in
H^0(M,L^N)$, we have
$$\var\big(Z_{s_N},\phi\big) = N^{-m}
\left[\kappa_m\,\|\ddbar\phi\|_2^2
 +O(N^{-\half +\ep})\right]\;,$$ where $$\kappa_m =
 \frac{\pi^{m-2}}{4}\;\zeta(m+2)\;.$$
\end{theo}

\noindent Here, $\|\ddbar\phi\|_2$ denotes the $\lcal^2$ norm of
$\ddbar\phi$, i.e.\ writing $i\ddbar\phi= \psi \frac 1{m!} \om^m$, we
have $\|\ddbar\phi\|_2^2= \int \psi^2\frac 1{m!} \om^m = \int i\psi
\ddbar\phi $. (Of course, we may assume that
$\phi$ is of bidegree
$(m-1,m-1)$, since ($Z_{s_N},\phi)=0$ for forms $\phi$ of other
bidegrees.)

In particular, for the case $\dim M=1$, we note that
$|\ddbar\phi| = \half |\Delta \phi|$, and thus
\begin{equation}\var\big(Z_{s_N},\phi\big) = N\inv
\left[\frac{\zeta(3)}{16\pi}\,\|\Delta\phi\|_2^2
 +O(N^{-\half +\ep})\right]\;.\label{smooth1}\end{equation} The
 leading term in \eqref{smooth1} was obtained by Sodin
and Tsirelson \cite{ST} for the case of random polynomials $s_N\in
H^0(\CP^1,\ocal(N))$ and random holomorphic functions on $\C$ and
on the disk. (The constant $\frac{\zeta(3)}{16\pi}$ was given in a
private communication from M. Sodin.)

Our final result is an asymptotic normality result of the type
proved in \cite{ST}. It follows very easily from the analysis
underlying Theorem \ref{sharp} together with a general asymptotic
normality result of Sodin-Tsirelson.

\begin{theo} \label{AN} With the same notation and hypotheses as
in Theorem \ref{sharp}, the distributions of the random variables
$$ \frac{(Z_{s_N},\phi)-\E(Z_{s_N},\phi)}{\sqrt{\var(Z_{s_N},\phi)}}$$
converge weakly to the
standard Gaussian distribution
$\ncal(0, 1)$ as $ N \to \infty$.\end{theo}

We let $\ncal(0, \sigma)$ denote the (real)
Gaussian distribution of mean zero and variance
$\sigma^2$. Substituting the values of the expectation and variance of
$(Z_{s_N},\phi)$ from \eqref{ave} and Theorem \ref{sharp},
respectively, we have
\begin{cor} With the same hypotheses as
in Theorem \ref{sharp}, the distributions of the random variables
$
N^{m/2}(Z_{s_N}-N\om,\phi)$ converge weakly to 
$\ncal(0, \sqrt{\kappa_m}\, \|\ddbar\phi\|_2)$ as $ N \to \infty.$
\end{cor}

Let us briefly mention some key ideas in the proofs and the
relation of the Sodin-Tsirelson methods to ours. The
Sodin-Tsirelson estimate  was  based on their construction of a
`bi-potential' for the pair correlation measures, i.e.\ functions
$ G_N(z, w)$ such that
\begin{equation}\label{biST}\Delta_z \Delta_w  G_N (z,w) =
K_2^N(z,w). \end{equation} Here, $K_2^N$ is the `pair correlation
function' for the zeros of sections of $\scal_N$, that is, the
probability density that a section  in $\scal_N$ has zeros at two
points $z$ and $w$ of $\CP^1$. The bi-potential is given in
\cite{ST} as a power series in the \szego kernel for
$\ocal(N)\to\CP^1$. (Here, the notation in (\ref{biST}) is taken
from \cite{BSZ1} and is not used in \cite{ST}.) In fact, the same
bipotential already  arose in \cite{BSZ1} in the setting of line
bundles over  a compact \kahler manifold as  the bi-potential for
the `pair correlation current,' i.e.
\begin{equation}\label{bi}-\d_z\dbar_z \d_w\dbar_w  G_N (z,w) =
\E\left( Z_{s_N}(z) \otimes  Z_{s_N}(w)\right)\;.\end{equation} We
build on our analysis of this bi-potential  in \cite{BSZ1} to
prove Theorem \ref{sharp}.

The other main ingredient in the proofs  are estimates derived
from  the off-diagonal asymptotics of the \szego kernel in
\cite{SZ2}.  For the sake of completeness, we review
the derivation of these off-diagonal asymptotics in \S \ref{off}.
Off-diagonal estimates of the \szego kernel with sharper
(exponentially small) remainder estimates are  given in \cite{DLM,
MaMa}, but the estimates of \cite{SZ2} already suffice for our
applications.

Although we are emphasizing positive line bundles over compact
\kahler manifolds, our results (and their proofs) extend with no
essential change to   positive line bundles over noncompact
\kahler manifolds for which the orthogonal projector $\Pi_N$ onto
$\lcal^2 H^0(M, L^N)$, the $L^2$ holomorphic sections of a positive
line bundle with respect to a Hermitian metric and the \kahler
volume form, has analytic properties similar to those in the
compact case. A model for a positive line bundle over a
noncompact \kahler manifold is provided by the Heisenberg  line
bundle $L_\H \to \C^m$ associated to the reduced Heisenberg group by
the identity character, as described in detail in  \cite{BSZ2,
BSZ0}. In this case, the analogue of Theorem \ref{volume} is an asymptotic formula
(Corollary \ref{dilates}) for the volume variance of the zeros of random Gaussian
entire functions on the dilates
$\sqrtn\,U$ of a domain $U\subset\C^m$ . Other model examples are
given by  homogeneous Hermitian line bundles over bounded symmetric domains with
curvature equal to the Bergman \kahler metric. We briefly discuss the extension to
random holomorphic sections in the noncompact case in
\S
\ref{NONCOMPACT}.

It will readily be recognized that the variance and normality
problems in higher dimensions make sense for the simultaneous
zeros of $k$ independent sections $s_1, \dots, s_k$  of a line
bundle over an $m$-dimensional complex manifold and are perhaps
most interesting for the full codimension case $k = m$. The same
problem may be posed for the critical points of a single Gaussian
random section. However, new technical ideas seem to be necessary
to obtain limit formula for the intersections of  the random zero
currents $Z_{s_j}$. We hope to return to this problem elsewhere.

In conclusion, we thank M. Sodin for  discussions of his work with
B. Tsirelson on number variance and asymptotic normality for
random analytic functions of one variable.

\section{Expected distribution of zeros and \szego kernels}\label{EDZ}

In this section, we review the basic formula from
\cite{BSZ1,BSZ2,SZ} for the expected distribution of zeros of
Gaussian  random sections of holomorphic line bundles. We state it
here in a  general framework which we shall use in our forthcoming
paper on zeros of random fewnomials \cite{SZ3}.

 We let  $(L,h)$ be a Hermitian holomorphic
line bundle over a complex manifold $M$ (not necessarily compact),
and  let $\scal$ be a finite-dimensional subspace of $H^0(M,L)$.
We suppose that $\dim \scal\ge 2$ and we give $\scal$ a Hermitian
inner product. The inner product induces the complex Gaussian
probability measure
\begin{equation}\label{gaussian}d\gamma(s)=\frac{1}{\pi^m}e^
{-|c|^2}dc\,,\qquad s=\sum_{j=1}^{n}c_jS_j\,,\end{equation} on
$\scal$, where $\{S_j\}$ is an orthonormal basis for $\scal$ and
$dc$ is $2n$-dimensional Lebesgue measure. This Gaussian is
characterized by the property that the $2n$ real variables $\Re
c_j, \Im c_j$ ($j=1,\dots,n$) are independent complex Gaussian random
variables with mean 0 and variance 1; i.e.,
$$\E c_j = 0,\quad \E c_j c_k = 0,\quad  \E c_j \bar c_k =
\de_{jk}\,.$$

 We let
\begin{equation}\label{sdef}\Pi_\scal(z,z) = \sum_{j=1}^n
\|S_j(z)\|_h^2\;,\qquad z\in M\;,\end{equation} denote the {\it
\szego kernel\/} for $\scal$ on the diagonal. (See
\S\ref{s-powers} for a discussion of the \szego kernel.) We now
consider a local holomorphic frame $e_L$ over a trivializing chart
$U$, and we write $S_j = f_j e_L$ over $U$. Any section
$s\in\scal$ may then be written as
$$s = \langle c, F \rangle e_L^{\otimes N}\;, \quad \mbox{where\ \ \ }
F=(f_1,\dots,f_k)\;,\quad\langle c,F \rangle = \sum_{j = 1}^n c_j
f_j\;.$$ If  $s = f e_L$,  its Hermitian norm is given by
$\|s(z)\|_h = a(z)^{-\half}|f(z)|$ where \begin{equation}
\label{a} a(z) = \|e_L(z)\|_h^{-2}\;. \end{equation} Recall that
the curvature form of $(L,h)$ is given locally by
$$\Theta_h= \ddbar \log a\;,$$ and the
{\it Chern form\/} $c_1(L,h)$ is given by
\begin{equation}\label{chern}c_1(L,h)=\frac{\sqrt{-1}}{2 \pi}
\Theta_h=\frac{\sqrt{-1}}{2 \pi}\d\dbar\log a\;.\end{equation} The
current of integration $Z_s$ over the zeros of $s = \langle c,F
\rangle e_L$ is then given locally by the {\it Poincar\'e-Lelong
formula\/}:
\begin{equation} Z_s =
\frac{\sqrt{-1}}{ \pi } \partial \bar{\partial}\log | \langle c,F
\rangle|\;. \label{Zs} \end{equation} It is of course independent
of the choice of local frame $e_L$ and  basis $\{S_j\}$.

We now state our formula for the expected zero divisor for the
linear system $\scal$:

\begin{prop}\label{EZ}Let $(L,h)$ be a Hermitian line bundle on a complex
manifold $M$, and let $\scal$ be a finite dimensional subspace of $
H^0(M,L)$.  We give $\scal$ a Hermitian inner
product and we let $\ga$ be the induced  Gaussian probability measure
on $\scal$. Then the expected zero current of a random section
$s\in\scal$ is given by

\begin{eqnarray*}\E_\ga(Z_s)  &=&\frac{\sqrt{-1}}{2\pi}
\partial
\bar{\partial} \log \Pi_{\scal}(z, z)+c_1(L,h)\;.\end{eqnarray*}
\end{prop}

\begin{proof} Let
$\{S_j\}$ be  an orthonormal basis of $\scal$. As above, we choose
a local nonvanishing section $e_L$ of $L$ over $U\subset M$, and
we write
$$s=\sum_{j=1}^{n}c_jS_j=\langle c,F\rangle e_L\;,$$ where
$S_j=f_j
 e_L,\ F=(f_1,\dots,f_{k})$.  As in the proof of \cite{SZ},
Proposition 3.1, we then write $F(x)= |F(x)| u(x)$ so that $|u|
\equiv 1$ and
\begin{equation}\label{2terms}\log  | \langle c, F \rangle| = \log |F| +
\log  | \langle c, u \rangle|\;.\end{equation}
A key point is that $\E\big(\log  | \langle c, u \rangle|\big)$
is independent of $z$ (and in fact, is a universal constant
depending only on $n$), and hence $\E\big(d\log  | \langle c, u
\rangle|\big)=0$.

Thus by (\ref{Zs}), we have
\begin{eqnarray*}\big(\E_\ga(Z_s),\phi\big)&=&\frac{\sqrt{-1}}{ \pi}
\int_{\C^n} \left(\log  | \langle c,F \rangle| , \ddbar\phi\right)
d\gamma(c)\\ &=&\frac{\sqrt{-1}}{ \pi} \int_{\C^n} \left(\log  |F| ,
\ddbar\phi\right) d\gamma(c) +\frac{\sqrt{-1}}{ \pi}\int_{\C^n}\left( 
\log  | \langle c, u
\rangle|, \ddbar\phi\right) d\gamma(c) \;,
\end{eqnarray*}
for all test forms $\phi\in\dcal^{m-1,m-1}(U)$. 
 The first term is independent of $c$
so we may remove the Gaussian integral. The vanishing of the second
term follows by noting that 
\begin{eqnarray*}\int_{\C^n}\left( 
\log  | \langle c, u
\rangle|, \ddbar\phi\right) d\gamma(c) &=& \int_{\C^n}d\gamma(c)
\int_M \log  | \langle c, u
\rangle|\, \ddbar\phi\\&=& \int_M\left[ \int_{\C^n}
 \log  | \langle c, u
\rangle|d\gamma(c)\right] \ddbar\phi \ =
\ 0\;,
\end{eqnarray*}
since $
\int
\log  | \langle c, u \rangle|\, d\gamma(c)= \frac 1\pi\int_\C
\log  | c_1|\,e^{-|c_1|^2} dc_1$ is constant, by the
$\U(n)$-invariance of $d\ga$. Fubini's Theorem can be applied above
since
$$\int_{M\times\C^n}\left|\log  | \langle c, u
\rangle|\, \ddbar\phi\right|\;d\gamma(c)= \left(\frac 1\pi\int_\C
\big|\log  | c_1|\big|\,e^{-|c_1|^2} dc_1\right)
\left(\int_M |\ddbar\phi|\right)< +\infty\;.$$

Thus
\begin{equation}\label{Egamma}\E_\ga(Z_s)=\frac{\sqrt{-1}}{2 \pi} \partial
\bar{\partial}\log  |F|^2  =\frac{\sqrt{-1}}{2 \pi} \partial
\bar{\partial}\left(\log \sum\|S_j\|_h^2+\log
a\right)\;.\end{equation} Recalling that $\Pi_{\scal}(z,
z)=\sum\|S_j(z)\|_h^2$ and that $c_1(L,h)=\frac{\sqrt{-1}}{2
\pi}\d\dbar\log a$, the formula of the proposition follows.
\end{proof}

\begin{rem} The complex manifold $M$, the line bundle $L$ and
space $\scal$ as well as its inner product in Proposition \ref{EZ}
are all completely arbitrary.  We do not assume that $M$ is
compact or that $(L,h)$ has positive curvature.  We do not even
assume that $\scal$ is base point free. If $\scal$ has no base points (points
where all sections in $\scal$ vanish), then we have the alternate
formula (see \cite{SZ})
$$ \E_\ga(Z_s)= \Phi^*_\scal\om_\FS\;,$$ where $\Phi_\scal:M\to
\PP\scal^*$ is the Kodaira map and $\om_\FS$ is the Fubini-Study
form on $\PP\scal^*$. In the general case where there are base
points, we have
$$ \E_\ga(Z_s)= \Phi^*_\scal\om_\FS +D\;,$$
where $D$ is the fixed component of the linear system $\PP\scal$.
\end{rem}

\subsection{Powers of an ample line bundle}\label{s-powers}

We now let $L\to M$ be an ample line bundle on a compact complex
manifold $M$.  We consider tensor powers $L^N=L^{\otimes N}$ of
the line bundle, and we let $\scal=H^0(M,L^N)$.  We further choose
a Hermitian metric $h$ on $L$ with strictly positive curvature and
we give $M$ the \kahler form $\om= \frac i2\Theta_h=\pi c_1(L,h)$.

We now describe the natural Gaussian probability measures on
$H^0(M,L^N)$ used in \cite{SZ,SZ2,BSZ1,BSZ2}.  For the case of
polynomials in one variable, these Gaussian ensembles are
equivalent to the $\SU(2)$ ensembles studied in \cite{BBL, Han,
NV, SZ} and elsewhere.

\begin{defin}\label{HG} Let $(L,h)\to (M,\om)$ be as above,
and let $h^N$ denote the Hermitian metric on $L^N$ induced by $h$.
We give $H^0(M,L^N)$ the inner product induced by the \kahler form
$\om$ and the Hermitian metric $h^N$:
\begin{equation}\label{inner}\langle s_1, \bar s_2 \rangle = \int_M h^N(s_1,
s_2)\,\frac 1{m!}\om^m\;,\qquad s_1, s_2 \in
H^0(M,L^N)\,\;.\end{equation}  The \ {\em Hermitian Gaussian measure}
on
$H^0(M,L^N)$ is the complex Gaussian probability measure $\ga_N$
induced by the inner product \eqref{inner}:
\begin{equation*}d\ga_N(s)=\frac{1}{\pi^m}e^
{-|c|^2}dc\,,\qquad s=\sum_{j=1}^{d_N}c_jS^N_j\,,\end{equation*}
where $\{S_1^N,\dots,S_{d_N}^N\}$ is an orthonormal basis for
$H^0(M,L^N)$.
\end{defin}

As in \cite{SZ, BSZ1} and elsewhere,  we analyze  the \szego
kernel for $H^0(M,L^N)$ by lifting it to the circle bundle
$X{\buildrel {\pi}\over \to} M$  of unit vectors in the dual
bundle $L\inv\to M$ with respect to $h$. In the standard way (loc.
cit.), sections of $L^N$ lift to equivariant functions on $X$.
Then $s\in H^0(M,L^N)$ lifts to a $CR$ holomorphic functions on
$X$ satisfying  $\hat s(e^{i\theta}x)= e^{iN\theta}\hat s(x)$. We
denote the space of such functions by  $ \hcal^2_N(X)$. The {\it
\szego projector\/} is the orthogonal projector
$\Pi_N:\lcal^2(X)\to\hcal^2_N(X)$, which is given by the {\it
\szego kernel}
$$\Pi_N(x,y)=\sum_{j=1}^{d_N} \wh S^N_j(x)\overline{\wh S^N_j(y)}
\qquad (x,y\in X)\;.$$ (Here, the functions $\wh S^N_j$ are the
lifts to $\hcal^2_N(X)$ of the orthonormal sections $S_j^N$; they
provide an orthonormal basis for $\hcal^2_N(X)$.)

Further, the covariant derivative $\nabla s$  of a section $s$ lifts
to the horizontal derivative $\nabla_h \hat{s}$ of its
equivariant lift $\hat{s}$ to $X$; the horizontal derivative is of the
form
\begin{equation}\label{HORDER} \nabla_h \hat s
=\sum_{j=1}^m\left( \frac{\d \hat s}{\d z_j} -A_j\frac{\d\hat s}{\d
\theta}\right)dz_j.
\end{equation}
For further discussion and  details on lifting sections, we refer
to \cite{SZ}.

We shall write $$|\Pi_N(z,w)|:=|\Pi_N(x,y)|\;,\quad
z=\pi(x),\,w=\pi(y)\in M\;.$$ In particular, on the diagonal we
have $\Pi_N(z,z)=\Pi_N(x,x)$, where $\pi(x)=z$. Note that
$\Pi_N(z,z)=\Pi_\scal(z,z)$ as defined by (\ref{sdef}) with
$\scal=H^0(M,L^N)$.

\subsection{\label{NONCOMPACT} Random functions and \szego kernels on
noncompact domains}

As mentioned in the introduction, Theorems
\ref{number}--\ref{AN}  extend with no
essential change to positive line bundles over noncompact
complete \kahler manifolds as long as the orthogonal projection
onto the space $\lcal^2H^0(M, L^N)$ of $L^2$ holomorphic sections  with
respect to the inner product (\ref{inner}) possesses the
analytical properties stated in  Theorem \ref{near-far} (and
mostly proved in \cite{SZ2}) for \szego kernels in the compact
case. It would take us too far afield to discuss in detail the
properties of \szego kernels and random holomorphic sections in
the noncompact setting, but we can illustrate the ideas with
homogeneous models.

Before discussing our specific noncompact models, we first note that
Proposition \ref{EZ} holds for infinite-dimensional spaces of
Gaussian random holomorphic sections.  There are several equivalent
ways to describe Gaussian random analytic sections or functions in an
infinite dimensional space (e.g., \cite{J, GJ,
ST}). To take a simple approach, we suppose that
$\{S_1,S_2,\dots,S_n,\dots\}$ is an infinite sequence of holomorphic
sections of a Hermitian line bundle $(L,h)$ on a (noncompact) complex
manifold $M$ such that 
\begin{equation}\label{unifL2} \sup_{z\in K}\sum_{j=1}^\infty
\|S_j(z)\|_h^2 <+\infty \quad \mbox{for all compact }\ K\subset
M\;.\end{equation}  We then consider the ensemble $(\scal,d\ga)$ of
sections of $L$ of the form
\begin{equation}\label{infinite} \scal=\left\{s=\sum_{j=1}^\infty
c_jS_j:c_j\in\C\right\}\;,\qquad d\ga= \prod_{j=1}^\infty\left(
\frac 1\pi e^{-|c_j|^2} dc_j\right)\;,\end{equation} i.e.\ we consider
random sections
$s= \sum_{j=1}^\infty
c_jS_j$, where the
$c_j$ are i.i.d.\ standard complex Gaussian random variables. It is
well known that \eqref{unifL2} implies that the series in
\eqref{infinite} almost surely converges uniformly on compact sets (see
e.g.\
\cite{J, Kah}), and hence with probability one,  $s\in H^0(M,L)$. We
then have:

\begin{prop}\label{EZinfinite} The expected zero current of the random
section
$s\in\scal$ in \eqref{infinite} is given by

$$\E(Z_s)  =\frac{\sqrt{-1}}{2\pi}
\partial
\bar{\partial} \log \Pi(z, z)+c_1(L,h)\;,$$
where $$\Pi(z, z)=\sum_{j=1}^\infty
\|S_j(z)\|_h^2\;.$$
\end{prop}
\begin{proof} The proof is exactly the same as the proof of
Proposition \ref{EZ} using the ensemble \eqref{infinite} with the
infinite product measure,  except we cannot use unitary invariance to
show that 
\begin{equation}\label{invariance} \textstyle
\int\log  | \langle c, u \rangle|\, d\gamma(c)= \frac 1\pi\int_\C
\log  | c_1|\,e^{-|c_1|^2} dc_1,\ \int
\big|\log  | \langle c, u \rangle|\big| d\gamma(c)= \frac 1\pi\int_\C
\big|\log  | c_1|\big|\,e^{-|c_1|^2} dc_1.\end{equation} To verify
 \eqref{invariance} in this case, we note
that $ \langle c, u(z) \rangle$ is a complex Gaussian random variable
of mean 0 and variance 1 (see \cite{J,Kah}), and hence   
$$\int f( \langle c, u \rangle)\, d\gamma(c)= \frac 1\pi\int_\C
f(\zeta)\,e^{-|\zeta|^2}\,d\zeta\;, \quad \mbox{for all }\ 
f\in\lcal^1(\C,e^{-|\zeta|^2}d\zeta)\;.$$ The identities of
\eqref{invariance} then follow by letting  
$f(\zeta)= \log |\zeta|$, resp. $f(\zeta)= \big|\log |\zeta|\big|$. 
\end{proof}

We are interested in the case where $(L,h)$ has positive curvature,
$M$ is complete with respect to the \kahler metric $\om=\frac i2
\Theta_h$, and
$\{S_j\}$ is an orthonormal basis of $\lcal^2H^0(M,L)$ with respect to 
the inner product \eqref{inner}.  Note that with probability one, a
random section $s$ is not an  $\lcal^2$ section (since
$\|s\|_2=\|c\|_2=+\infty\ a.s.$), but is a holomorphic section of
$L$. (Equivalently,
$\lcal^2 H^0(M,L)$ carries a Gaussian measure in the sense of
Bochner-Minlos; see \cite{GJ}.)

The first model noncompact  case is known as the
Bargmann-Fock space
$$\fcal:= \hcal^2(\C^m, e^{-|z|^2}) =\left \{f\in\ocal(\C^m):
\int_{\C^m}|f|^2e^{-|z|^2}\,dz <+\infty\right\}.$$ We can regard
elements of $\fcal$ as $\lcal^2$ sections of the trivial bundle
$L_\H$ over $\C^m$ with metric $h=e^{-|z|^2}$.  The associated
circle bundle $X$ can be identified with the reduced Heisenberg
group; see \cite[\S 2.3]{BSZ0} or \cite[\S 1.3.2]{BSZ2}.  Then
$\fcal= \lcal^2 H^0(\C^m,L_\H)=\hcal^2_1(X)$, and more generally,
\begin{equation} \label{Heis} \lcal^2H^0(\C^m,L^N_\H) = \hcal_N^2(X)
=\left
\{f\in\ocal(\C^m): \int_{\C^m}|f|^2e^{-N|z|^2}\,dz
<+\infty\right\}\;.\end{equation} In dimension one, this example
is referred to as the `flat model'
 in \cite{ST}.

An orthonormal basis for the Hilbert space $\lcal^2H^0(\C^m,L^N_\H)$
with inner product $\langle f_1, \bar f_2\rangle = \int_{\C^m}f_1\bar
f_2^2e^{-|z|^2}\,dz$ is $$\left\{S^N_k(z)=\frac
{N^{m/2}}{\pi^{m/2}} \frac {(N^{1/2}z)^k}{\sqrt{k!}}
\right\}_{k\in \N^m},
$$ where we use the usual conventions $z^k=z_1^{k_1}\cdots z_m^{k_m},\
k!=k_1!\cdots k_m!,\ |k|=k_1+\cdots +k_m$.
A Gaussian random section is defined
by \begin{equation}\label{BFgauss}
f_N(z) = \sum_{k \in \N^m} c_k S^N_k(z) =\frac
{N^{m/2}}{\pi^{m/2}}  \sum_{k \in \N^m}\frac{c_k}{\sqrt{k!}}
(N^{1/2}z)^k\;,\end{equation} where the coefficients $c_k$ are
independent  standard complex Gaussian random variables as in
(\ref{gaussian}). As mentioned above, the random sections $f_N$ are
almost surely not in
$\lcal^2H^0(\C^m,L^N_\H)$.  However, they are almost surely entire
functions of finite order 2 in the sense of Nevanlinna theory.  Indeed,
we easily see from \eqref{infinite} that
$$\ga\left(\left\{c\in\C^\infty: |c_k|^2 \le \textstyle2\sum_j\log
k_j\ \ \mbox{for }\ k_j\ge 2\right\}\right) >0 \;,$$ and hence it
follows from
 the zero-one law that
$|c_k|^2 =O\left(\sum_j
\log k_j\right)\ a.s.$ Therefore, by Cauchy-Schwartz, $$|f(z)| =O
\left(\sum (1-\ep)^{|k|}\right)^{1/2}  \left(\sum \frac
{(1+2\ep)^{|k|}}{k!}N^{|k|}|z^k|^2\right)^{1/2} =
O\left(e^{(\half+\ep)N|z|^2}\right) \quad a.s.$$ for all $\ep>0$. 
Thus we have an upper bound for the Nevanlinna growth function, 
$$T(f_N,r):= \mbox{Ave}_{\{|z|=r\}}\log^+|f(z)| \le \sup_{|z|\le
r}\log^+|f(z)|\le 
\left[\frac N2 +o(1)\right]r^2 \quad a.s.$$ (where $o(1)$ denotes a
term that goes to 0 as $r\to\infty$, for  each fixed $N\ge 1$). On the
other hand, if
$T(f_N,r) = O(e^{(\half-\ep)N|z|^2})$, then
$f_N\in\lcal^2H^0(\C^m,L^N_\H)$, which has probability zero.  Thus,
$$\limsup_{r\to\infty}\frac{T(f_N,r)}{r^2} = \frac N2 \qquad a.s.$$

To use the proofs in \S\S \ref{s-sharp}--\ref{s-normality} to show that
Theorems
\ref{number}--\ref{AN} hold for the line bundle
$L_\H$, we need only to verify that the \szego kernel $\Pi_N^\H$, i.e.\
the kernel of the orthogonal projection to $\lcal^2 H^0(\C^m,L_\H)$,
satisfies the diagonal and off-diagonal asymptotics in Theorem
\ref{near-far}. In the model Heisenberg case, the \szego kernel is
given by
\begin{equation}\label{heisen-N}\Pi^\H_N(z,\theta;w,\phi) =
e^{iN(\theta-\phi)}\sum_{k\in \N^m}S_k(z)\overline{S_k(w)} =
\frac{N^m}{\pi^m} e^{i  N(\theta - \phi)+N z \cdot \bar{w}- \frac
N2 (|z|^2+ |w|^2)} \;,\end{equation} (see \cite{BSZ2}) and visibly has
these properties. 

Another class of homogeneous examples are the bounded symmetric
domains $\Omega \subset \C^m$, equipped with their Bergman metrics
$\omega =\frac i2 \ddbar \log K(z, \bar{z})$ where $K(z, \bar{z})$
denotes the Bergman kernel function of $\Omega$. Let $(L, h) \to \Omega$
be the holomorphic homogeneous Hermitian line bundle over $\Omega$
with curvature $(1,1)$ form $\omega$. It was observed by Berezin
\cite{Ber} that the \szego kernels $\Pi_N$ for $\lcal^2H^0(\Omega,
L^N)$ also have the form $C_N e^{N \psi(z, \bar{w})}$ where $C_N$
is a normalizing constant and $\psi = \log K(z, \bar{w})$. In the
case of the unit disc $D \subset \C$ with its Bergman (hyperbolic
metric) $\frac{-i}{2}\ddbar \log (1 - |z|^2)$, the space $\lcal^2
H^0(D, L^N)$ may be identified with the holomorphic discrete
series irreducible representation $\dcal_N^+$ of $SU(1,1)$ (cf.\
\cite[p.~40]{K}), that is with the space of holomorphic functions
on $D$ with inner product
$$\|f\|_N^2 = \int_{D} |f(z)|^2 (1 - |z|^2)^{N-2} dz. $$
The  factor $e^{N \log (1 - |z|^2)}$ comes from the Hermitian
metric.
 An
 orthonormal
basis for the holomorphic sections of $L^N$ is then given by the
monomials ${N + n - 1 \choose n}^{1/2} z^n\;$ ($n = 0, 1, 2,
\dots$). The \szego kernels are given by  $\Pi_N(z,w)= (1 - z
\bar{w})^N$.  The \szego kernels also visibly have the properties
stated in Theorem \ref{near-far}. These ensembles are called the
hyperbolic model in \cite{ST}. Random SU(1,1) polynomials are
studied in \cite{BR}, where further details can be found.

Thus our proofs also yield the
following result:
\begin{theo} Theorems \ref{number}--\ref{AN} hold for the zeros of 
sections in the following ensembles:
\begin{itemize} \item  Gaussian random sections $f_N\in
H^0(\C^m,L^N_\H)$ given by
\eqref{BFgauss};\item Gaussian random sections of the holomorphic
homogeneous Hermitian line bundle $(L,h)$ over a bounded
symmetric domain $ \Omega\subset\C^m$, as described above.
\end{itemize}
\end{theo}

We note that taking the $N$-th power of the line bundle $L_\H\to\C^m$
(i.e., taking the $N$-th power of the metric $e^{-|z|^2}$) corresponds
to dilating
$\C^m$ by $\sqrtn$.  Precisely, the map $$\tau_N:
\lcal^2H^0(\C^m,L_\H)\to \lcal^2H^0(\C^m,L^N_\H)\;,\qquad
(\tau_Nf)(z):= N^{m/2}f(N^{1/2}z)\;,$$ is unitary. Thus we can
restate our result on the volume (or number, in dimension 1) variance
for the Bargmann-Fock ensemble as follows:
\begin{cor}\label{dilates} Let
  $$f(z) =  \sum_{k \in \N^m}\frac{c_k}{\sqrt{k!}}\,
z^k\;,$$ where the coefficients $c_k$ are independent  complex
Gaussian random variables with mean 0 and variance 1.   Let $U$ be a
domain in
$\C^m$ with piecewise
$\ccal^2$ boundary and no cusps, and consider its dilates
$U_N:=\sqrtn\,U$.  Then,
$$\var\big(\vol_{2m-2}[Z_{f}\cap U_N]\big) =
\nu_m\,\vol_{2m-1}(\d U_N)
 +O(N^{-\half +\ep})\;,$$ where $\nu_m =
\frac{\pi^{m-5/2}}{8}\;\zeta(m+\textstyle \half)$.
\end{cor}
\noindent Note that $\d U_N =N^{m-1/2} \d U$, so we have
$\var\big(\vol_{2m-2}[Z_{f}\cap \sqrtn\,U]\big)\sim N^{m-1/2}$.

Off-diagonal estimates for general Bergman kernels of positive
line bundles over complete \kahler manifolds  are proved in
\cite{MaMa} using  heat kernel methods.
  The  relevant  issue for
 this article is the approximation of the $L^2$ \szego
kernel by its Boutet de Monvel -Sj\"ostrand parametrix in the
noncompact case. The analysis of \szego kernels on noncompact
spaces  lies outside the scope of this article, so we do not state
the general results here. But it appears that the general results
of \cite{MaMa}  give sufficient control over \szego kernels in the
noncompact case to allow Theorems
\ref{number}--\ref{AN} to be extended to
all positive line bundles over complete \kahler manifolds.

\section{A bipotential for the variance}

Our proofs  of Theorems \ref{number}-- \ref{sharp} are based on a
bipotential implicitly given in \cite{SZ}. To describe our
bipotential
$Q_N(z,w)$, we define the function
\begin{equation}\label{Gtilde} \wt G(t):= -\frac 1{4\pi^2}
\int_0^{t^2} \frac{\log(1-s)}{s}\,ds\;, \qquad 0\le t\le
1.\end{equation}  Alternately,
\begin{equation}\label{Gtilde1}\wt G(e^{-\la}) = -\frac
1{2\pi^2} \int_\la^\infty
\log(1-e^{-2s})\,ds\;,\qquad \la>0\;.\end{equation}
(The function  $\wt G$ is a modification of the function
$G$ defined in \cite{SZ}; see \eqref{GGtilde}.) We also
introduce the {\it normalized \szego kernel}
\begin{equation}\label{PN} P_N(z,w):=
\frac{|\Pi_N(z,w)|}{\Pi_N(z,z)^\half
\Pi_N(w,w)^\half}\;.\end{equation}

\begin{defin} Let $(L,h)\to(M,\om)$ be as in Theorems
\ref{volume}--\ref{sharp}.  The \ {\em variance bipotential} is the
function
$Q_N:M\times M\to [0,+\infty)$ given by
\begin{equation}
\label{QN} Q_N(z,w)= \wt G(P_N(z,w)) = -\frac 1{4\pi^2}
\int_0^{P_N(z,w)^2} \frac{\log(1-s)}{s}\,ds\;. \end{equation}
\end{defin}

We remark here that $Q_N$ is $\ccal^\infty$ off the diagonal, but is only
$\ccal^1$ and not $\ccal^2$ at all points on the diagonal in $M\times M$, as
the computations in \S \ref{s-number} show.

The variance in Theorems \ref{number}--\ref{sharp} can be given by a double
integral of the bipotential, as stated in Propositions \ref{BIPOT} and
\ref{varint2} below:

\begin{prop} \label{BIPOT} Let $(L,h)\to (M,\om)$ and
$\phi$ be a $(2m-2)$-form on $M$ with $\ccal^2$ coefficients. Then
\begin{equation} \label{varint1} \var( Z_{s_N},\phi) =
\int_M \int_M Q_N(z,w)\,(i\ddbar \phi(z))\, (i\ddbar \phi(w))
\;.\end{equation}
\end{prop}

To begin the proof of the proposition, we write
\begin{equation}\label{Psi}\Psi_N=(S_1^N,\dots,S_{d_N}^N)\in
H^0(M,L^N)^{d_N}\;,\end{equation} where $\{S_j^N\}$ is an
orthonormal basis of $H^0(M,L^N)$. As in the proof of Proposition
\ref{EZ}, we write $$\Psi_N(z) = |\Psi_N(z)|\, u_N(z)\;,$$ where
$|\Psi_N|:= (\sum_j \|S^N_j\|_{h^N}^2)^{1/2}$, so that
$|u_N|\equiv 1$.

We first establish a less explicit variance formula:

\begin{lem} \label{varinta}
$$\var( Z_{s_N},\phi) = \frac{1}{\pi^2}
 \int_M \int_M (i\ddbar \phi(z)) (i\ddbar \phi(w))
\int_{\C^{d_N}} \log |\langle u_N(z), c\rangle| \,\log |\langle
u_N(w), c\rangle|\, d\ga_N(c)\;.$$
\end{lem}

\begin{proof}

 We write sections $s_N\in
H^0(M,L^N)$ as
\begin{equation}\label{Psi0}s_N=\sum_{j=1}^{d_N} c_jS_j^N = \langle
c, \Psi_N\rangle\;,\qquad c=(c_1,\dots,c_{d_N}) \;.\end{equation}
Writing $\Psi_N=Fe_L^{\otimes N}$, where $e_L$ is a local
nonvanishing section of $L$, and recalling that $$\om=\frac
i2\Theta_h=- i \ddbar\log \|e_L\|_h\;,$$ we have by \eqref{Zs},
\begin{eqnarray}\label{ZsN}Z_{s_N}&=&\frac i\pi \ddbar \log
|\langle c,F\rangle| \ =\ \frac i\pi \ddbar \log |\langle
c,\Psi_N\rangle| -\frac i\pi \ddbar \log\|e_L^{\otimes
N}\|_h\nonumber\\& =& \frac i\pi \ddbar \log |\langle
c,\Psi_N\rangle|+ \frac N\pi \,\om\;.\end{eqnarray}

Let $\phi$ be a test form, and consider the
random variable
\begin{equation}\label{YN}Y_N:=\left( \frac i\pi \ddbar \log
|\langle c,\Psi_N\rangle|\,,\,\phi\right) = \left( \log
|\langle c,\Psi_N\rangle|\,,\,\frac i\pi \ddbar
\phi\right)\;,\end{equation} so that $(Z_{s_N},\phi) = Y_N
+NC$, where $ C=
\frac 1\pi\int_M
\om\wedge\phi$; hence $$\var (Z_{s_N},\phi)=\var (Y_N)\;.$$

By \eqref{Egamma}, we have
\begin{equation}\label{EYN} \E(Y_N) = \left( \frac i\pi \ddbar \log
|\Psi_N|\,,\,\phi\right) = \frac i\pi
\int_M\log|\Psi_N|\,\ddbar\phi \;,\end{equation} whereas by
\eqref{YN}, we have
\begin{equation}\label{EZV} \E( Y_N^2)= \frac{1}{\pi^2} \int_M
\int_M (i\ddbar \phi(z))\;(i\ddbar \phi(w)) \int_{\C^{d_N}} \log
|\langle c,\Psi_N(z)\rangle|\, \log |\langle c,\Psi_N(w)
\rangle|\, d\ga_N(c)\;.
\end{equation}
Recalling that $\Psi_N = |\Psi_N| u_N$ with $|u_N|\equiv 1$, we
have
\begin{eqnarray}\log |\langle \Psi_N(z), c\rangle| \,\log |\langle
\Psi_N(w), c\rangle| &=&
\log |\Psi_N(z)| \,\log |\Psi_N(w)| + \log|\Psi_N(z)| \,\log |\langle
u_N(w), c\rangle|\nonumber \\&&+ \log |\Psi_N(w)| \,\log |\langle u_N(z),
c\rangle |\nonumber \nonumber\\&&+
 \log |\langle u_N(w), c\rangle |  \,\log |\langle u_N(z), c\rangle
|\;,\label{4terms}\end{eqnarray} which decomposes (\ref{EZV}) into
four terms. By (\ref{EYN}), the first term contributes
\begin{equation}\frac{1}{\pi^2}  \int_M \int_M (i\ddbar \phi(z))
(i\ddbar \phi(w))\log |\Psi_N(z)| \,\log |\Psi_N(w)|  =  (\E
Y_N )^2\;.
\end{equation} The $c$-integral in the second term is independent of $w$ and hence
the second term vanishes.  The third term likewise vanishes.
Therefore, the fourth term  gives the variance $\var(
Z_{s_N},\phi)$.
\end{proof}

We now complete the proof of Proposition \ref{BIPOT} by evaluating
the $c$-integral of Lemma \ref{varinta}:

\begin{lem} \label{varintb} We have:
$$\frac 1{\pi^2}\int_{\C^{d_N}}\log |\langle u_N(z), c\rangle| \,
\log |\langle u_N(w), c\rangle|\, d\ga_N(c) = Q_N(z,w)+K\;,$$
where $K$ is a universal constant.
\end{lem}

\begin{proof}

We showed in  \cite[p.~779]{SZ} by an elementary computation that
\begin{equation}\label{equalsG}\int_{\C^{d_N}} \log |\langle u_N(z),
c\rangle| \,\log |\langle u_N(w), c\rangle|\, d\ga_N(c)=
G(|\langle u_N(z),\overline{ u_N(w)}\rangle|)\;,\end{equation}
where
\begin{equation}\label{G} G(t)= \frac{1}{\pi^2} \int_{\C^2
}e^{-(|c_1|^2 + |c_2|^2)}  \log |c_{1}|\,\log \left|c_1t + c_2
\sqrt{1-t^2} \right|\, dc_1 dc_2\;.\end{equation}

The computation of the integral \eqref{G} was begun in
\cite[\S4]{SZ}.  Let us finish it.  By (47)--(50) in \cite{SZ}
(with $\la=\half r^2$), we have
$$G(e^{-\la})=k_1+k_2 \la +\half
\int_0^\la \log(1-e^{-2s})\,ds\;.$$ Since $G(0)$ is finite,
$k_2=0$ and hence

\begin{equation}\label{finishG}
G(e^{-\la})=k_0-\half \int_\la^\infty \log(1-e^{-2s})\,ds\;,
\end{equation} or equivalently,
\begin{equation}\label{finishG1}
G(t)=k_0 -\frac 14 \int_0^{t^2} \frac{\log(1-s)}{s}\,ds \qquad (0\le t\le
1)\;.
\end{equation}
Hence,
\begin{equation}\label{GGtilde} \wt G(t)= \frac 1{\pi^2}
[G(t)-k_0]\;.\end{equation}

The lemma follows from \eqref{equalsG} and \eqref{GGtilde} with
\begin{equation}\label{t=PN}t = P_N(z, w)= \frac{|\langle \Psi_N(z),
\overline{\Psi_N(w)}\rangle|}{|\Psi_N(z)|\,|\Psi_N(w)|} = |\langle
u_N(z),\overline{ u_N(w)}\rangle|\;.\end{equation}
\end{proof}

Proposition \ref{BIPOT} is an immediate consequence of Lemmas
\ref{varinta}--\ref{varintb}.\qed

\medskip
To state the bipotential formula for the volume variance of Theorem
\ref{volume} (and number variance of Theorem \ref{number}), we let
\begin{equation}\label{Phi}\Phi:=
\frac 1 {(m-1)!}\om^{m-1}\;,\end{equation} so that
$\vol_{2m-2}[Z_{s_N}\cap U] = \big(Z_{s_N},\chi_U\,\Phi)$.  (If
$m=1$, we set $\Phi=1$.)

\begin{prop}\label{varint2} Let $(L,h)\to (M,\om),\ U\subset M$
be as in Theorem \ref{sharp}. Then \begin{equation*}
 \var\big(\vol_{2m-2}[Z_{s_N}\cap U]\big) =- \int_{\d U\times
\d U}\dbar_z \dbar_w Q_N(z,w)
\wedge\Phi(z)\wedge\Phi(w)\;,\end{equation*}  where
$Q_N$ is given by
\eqref{QN}.  In particular, for the case where  $\dim M=1$, we have
$$\var\big(\#[Z_{s_N}\cap U]\big) =- \int_{\d U\times
\d U}\dbar_z\dbar_wQ_N(z,w)\;.$$
 \end{prop}

\begin{proof}  As
in the proof of Lemma
\ref{varinta}, we let \begin{equation}\label{YNUm}Y_N^U:=\left(
\frac i\pi \ddbar \log |\langle
c,\Psi_N\rangle|\,,\,\chi_U\,\Phi\right) = \int _U\frac
i\pi\ddbar
\log |\langle c,\Psi_N\rangle|\wedge\Phi= \frac
i{\pi}\int_{\d U}
 \dbar\log |\langle c,\Psi_N\rangle|
\wedge\Phi\;,\end{equation} where
$\Psi_N$ is given by \eqref{Psi}. As
before, we have $(Z_{s_N},\chi_U\,\Phi)=  Y_N^U +NC$, where
$C=\frac 1\pi\int_U \om\wedge\Phi=\frac m\pi \vol(U)$.
Hence,
\begin{equation}\label{=var}\var\big(\vol_{2m-2}[Z_{s_N}\cap
U]\big) =\var (Z_{s_N},\chi_U\,\Phi)=\var
(Y_N^U)\;.\end{equation}

By \eqref{YNUm} (see also the proof of Proposition \ref{EZ}), we
then have
\begin{equation}\label{EYNU}
\E\big(Y^U_N\big)=\frac
i{\pi}\int_{\d U} \dbar\log |\Psi_N|
\wedge\Phi\;.\end{equation} Also by
\eqref{YNUm}, we have \begin{equation}\label{EZU}
\E\big((Y^U_N)^2\big)=
\frac{-1}{\pi^2} \int_{\C^{d_N}}\int_{\d U} \int_{\d U}
 \left[\dbar\log |\langle
c,\Psi_N(z)\rangle|\wedge\Phi(z)\right] \left[ \dbar\log |\langle
c,\Psi_N(w) \rangle|\wedge\Phi(w)\right]\,
d\ga_N(c)\;.\end{equation}Again following the proof of Lemma
\ref{varinta},
 we  write $\Psi_N = |\Psi_N| u_N$ with
$|u_N|\equiv 1$, and use \eqref{4terms} to decompose \eqref{EZU}
into four terms. The first term contributes
\begin{equation}\frac{-1}{\pi^2} \int_{\d U} \int_{\d U}
\left[\dbar\log |\Psi_N(z)|\wedge\Phi(z)\right] \left[
\dbar\log |\Psi_N(w) |\wedge\Phi(w)\right] = \left (\E
Y_N^U\right )^2\;.
\end{equation} The second term vanishes since  $\int
\log  | \langle c, u_N(w) \rangle|\, d\gamma_N(c)$ is
independent  of $w$ and hence $\dbar$ kills it. Similarly, the
third term vanishes, since it contains $\int
\log  | \langle c, u_N(z) \rangle|\, d\gamma_N(c)$, which is
independent of $z$.
Therefore,
\begin{eqnarray*}
\var( Y_N^U) &=&
\frac{-1}{\pi^2} \int_{\C^{d_N}}\int_{\d U} \int_{\d U}
 \left[\dbar\log |\langle u_N(z), c\rangle|\wedge\Phi(z)\right]
\left[
\dbar\log |\langle u_N(w), c\rangle|\wedge\Phi(w)\right]
d\ga_N(c)\\ &=&\frac{-1}{\pi^2} \int_{\C^{d_N}}
\dbar_z\dbar_w  \left(\int_{\d U\times \d U} \log |\langle
u_N(z), c\rangle| \,\log |\langle u_N(w),
c\rangle|\,\Phi(z)\wedge\Phi(w)\right)
d\ga_N(c)\,.\end{eqnarray*}
The
formula of the proposition then follows from Lemma \ref{varintb}.
\end{proof}

\section{Off-diagonal asymptotics and estimates for the \szego
kernel}\label{off}
In this section, we
use the off-diagonal asymptotics for $\Pi_N(z,w)$ from \cite{SZ2} to
provide the off-diagonal estimates for the  normalized \szego kernel
$P_N(z,w)$ that we need for our variance formulas.  Our estimates are of two
types: (1) `near-diagonal' asymptotics (Propositions
\ref{better}--\ref{best}) for
$P_N(z,w)$ where the distance
$d(z,w)$ between
$z$ and
$w$ satisfies an upper bound $ d(z,w)\le b\left(\frac
{\log N}{N}\right)^{1/2}$ ($b\in\R^+$); (2) `far-off-diagonal' asymptotics
(Proposition \ref{DPdecay}) where  $ d(z,w)\ge b\left(\frac
{\log N}{N}\right)^{1/2}$.

As discussed in \S \ref{s-powers} (cf.  \cite{Z, SZ, SZ2}), we
obtain the asymptotics by identifying the line bundle \szego
kernel $\Pi_N$ with a scalar \szego kernel $\Pi_N(x, y)$ on the
unit circle bundle $X \subset L^{-1} \to M$ associated to the
Hermitian metric $h$. Given $z_0 \in M$,  we choose a neighborhood
$U$ of $z_0$, a local normal coordinate chart
$\rho:U,z_0\to\C^m,0$ centered at $z_0$, and a {\it preferred\/}
local frame  at $z_0$, which we defined in \cite{SZ2} to be a
local frame $e_L$ such that
\begin{equation}\label{preferred}\|e_L(z)\|_h=1-\half \|\rho(z)\|^2
+ \cdots\;.\end{equation} For $u=(u_1,\dots,u_m)\in \rho(U),\
\theta\in (-\pi,\pi)$, we let
\begin{equation}\tilde\rho(u_1,\dots,u_m,\theta)=\frac
{e^{i\theta}} {|e^*_L(\rho\inv(u))|_h}
 e^*_L(\rho\inv(u))\in
X\,,\label{coordinates}\end{equation} so that
$(u_1,\dots,u_m,\theta)\in\C^m\times \R$ give local coordinates on
$X$. As in \cite{SZ2}, we write
$$\Pi_N^{z_0}(u,\theta; v,\phi)
=\Pi_N(\tilde\rho(u,\theta),\tilde\rho(v,\phi))\;.$$ Note that
$\Pi_N^{z_0}$ depends on the choice of coordinates and frame; we
shall assume that we are given normal coordinates and local frames
for each point $z_0\in M$ and that these normal coordinates and
local frames are smooth functions of $z_0$.

The scaling asymptotics of $\Pi_N^{z_0}(u,\theta; v,\phi)$ lead to
the model Heisenberg \szego kernel (\ref{heisen-N}) discussed in
\S \ref{NONCOMPACT}  for the Bargmann-Fock  space of functions on
$\C^m$. We shall use the following (near and far) off-diagonal
asymptotics from \cite {SZ2}:

\begin{theo} \label{near-far} Let $(L,h)\to (M,\om)$ be as in Theorem
\ref{sharp}, and let $z_0\in M$.  Then using the above notation,
\begin{enumerate}
\item[i)] $\  N^{-m}\Pi_N^{z_0}(\frac{u}{\sqrtn},\frac{\theta}{N};
\frac{v}{\sqrtn},\frac \phi N)$
$$\begin{array}{l}
= \Pi^\H_1(u,\theta;v,\phi)\left[1+ \sum_{r = 1}^{k} N^{-r/2}
p_{r}(u,v) + N^{-(k +1)/2} R_{Nk}(u,v)\right]\;,\end{array}$$
where
 the $p_r$ are polynomials in $(u,v)$ of degree
$\le 5r$ (of the same parity as
$r$), and
$$|\nabla^jR_{Nk}(u,v)|\le C_{jk\ep b}N^{\ep}\quad \mbox{for }\
|u|+|v|<b\sqrt{\log N}\,,$$ for
$\ep,b\in\R^+$,  $j,k\ge 0$. Furthermore, the constant  $C_{jk\ep
b}$ can be chosen independently of $z_0$. \item[ii)] For
$b>\sqrt{j+2k+2m}\,$, $j,k\ge 0$, we have
$$ \left|\nabla^j_h
\Pi_N(z,w)\right|=O(N^{-k})\qquad \mbox{uniformly for }\ d(z,w)\ge
b\,\sqrt{\frac {\log N}{N}} \;.$$
\end{enumerate}\end{theo}

Here $\nabla^j_h=(\nabla_h)^j$ is the $j$-th iterated horizontal
covariant derivative; see (\ref{HORDER}).  Theorem \ref{near-far} is
equivalent to equations (95)--(96) in \cite{SZ2}, where the result was
shown to hold for almost-complex symplectic manifolds. (The remainder
in (i) was given for $v=0$, but the proof holds without any change for 
$v\ne 0$.  Also the statement of the result was divided into the two
cases where the scaled distance is less or more, respectively,
than $N^{1/6}$ instead of $\sqrt{\log N}$ in the above
formulation, which is more useful for our purposes.)  A
description of the polynomials $p_r$ in part (i) is given in
\cite{SZ2}, but we only need the $k=0$ case in this paper. For the
benefit of the reader, we give a complete proof of Theorem
\ref{near-far} in \S\ref{s-proof} below.

\begin{rem} The \szego kernel actually satisfies the sharper `Agmon decay
estimate' away from the diagonal:
\begin{equation}\label{agmon}\nabla^j\Pi_N(
z,\theta;w,\phi)= O\left(e^{-A_j\sqrtn\,d(z,w)} \right)\;,
\qquad j\ge 0\;.
\end{equation} In particular,
\begin{equation}\label{szegodecay} |\Pi_N( z,w)|=  O\left(e^{-A\sqrtn\,d(z,w)}
\right)\;.\end{equation} A short proof of \eqref{szegodecay} is
given in \cite[Th.~2.5]{Be}; similar estimates were established by
M.~Christ \cite{Ch}, H. Delin \cite{D}, and N. Lindholm\cite{Li}.
(See also \cite{DLM, MaMa} for off-diagonal exponential estimates
in a more general setting.)  We do not need Agmon estimates for
this paper; instead Theorem \ref{near-far} suffices.
\end{rem}

We now state our far-off-diagonal decay estimate for $P_N(z,w)$, which
follows immediately from Theorem
\ref{near-far}(ii) and the fact that $\Pi_N(z,z)= \frac
1{\pi^m}N^m(1+O(N\inv))$ (by \cite{Z} or Theorem \ref{near-far}(i)).

\begin{prop}\label{DPdecay} Let $(L,h)\to (M,\om)$ be as in Theorem
\ref{sharp}, and let $P_N(z,w)$ be the normalized \szego kernel for
$H^0(M,L^N)$ given by \eqref{PN}. For
$b>\sqrt{j+2k}$,
$j,k\ge 0$, we have
$$ \nabla^j
P_N(z,w)=O(N^{-k})\qquad \mbox{uniformly for }\ d(z,w)\ge
b\,\sqrt{\frac {\log N}{N}} \;.$$
\end{prop}
\medskip
The normalized \szego kernel $P_N$ also satisfies Gaussian
decay estimates valid very close to the diagonal. To give the estimate, we write by
abuse of notation,
$$P_N(z_0+u,z_0+v):=
P_N(\rho\inv(u),\rho\inv(v))= \frac {|\Pi_N^{z_0}(u,0;
v,0)|}{\Pi_N^{z_0}(u,0; u,0)^{1/2}\Pi_N^{z_0}(v,0;
v,0)^{1/2}}\;.$$ As an immediate consequence of Theorem
\ref{near-far}(i), we have:

\begin{prop} \label{better} Let $P_N(z,w)$ be as in Proposition
\ref{DPdecay}, and let $ z_0\in M$.  For $b,\ep>0,\  j\ge 0$, there is a constant
$C_j=C_j({\ep,b})$, independent of the point $z_0$, such that
\begin{eqnarray*}\textstyle  P_N\left(z_0+\frac u{\sqrtn},z_0+\frac v{\sqrtn}\right) &=&
e^{-\half |u-v|^2}[1 + R_N(u,v)]\\ &&\quad |\nabla^jR_N(u,v)|\le C_j\,N^{-1/2+\ep}\quad
\mbox{for }\ |u|+|v|<b\sqrt{\log N}\;.\end{eqnarray*}
\end{prop}

 As a corollary we have:

\begin{prop} \label{best}
The remainder $R_N$ in Proposition \ref{better} satisfies
$$ |R_N(u,v)|\le \frac {C_2}2\,|u-v|^2N^{-1/2+\ep}, \quad |\nabla R_N(u)| \le
C_2\,|u-v|\,N^{-1/2+\ep},
 \quad \mbox{for }\ |u|+|v|<b\sqrt{\log N}.$$\end{prop}

\begin{proof} Since $P_N\left(z_0+u,z_0+v\right) \le 1=
P_N\left(z_0+u,z_0+u\right)$, we conclude that $R_N(u,u)=0,\ dR_N|_{(u,u)}=0$,
and thus by  Proposition
\ref{better},
$$|\nabla R_N(u,v)| \le  \sup_{0\le t\le
1}|\nabla^2 R_N(u,(1-t)u+tv)|\,|u-v|\le C_2\,|u-v|\,N^{-1/2+\ep}\;.$$
 Similarly,
$$|R_N(u,v)| \le \half  \sup_{0\le t\le
1}|\nabla^2 R_N(u,(1-t)u+tv)|\,|u-v|^2\le \frac {C_2}2\,|u-v|^2\,N^{-1/2+\ep}\;.$$
\end{proof}

\subsection{Proof of Theorem \ref{near-far}}\label{s-proof} In this section, we
sketch the proof of Theorem \ref{near-far}. The argument is
essentially contained in \cite{SZ2}, but we add some details
relevant to the estimates in Theorem \ref{near-far}.

The \szego kernels $\Pi_N(x,y)$ are the Fourier coefficients of
the total
\szego projector $\Pi(x,y):\lcal^2(X)\to \hcal^2(X)$; i.e.
$\Pi_N(x,y)=\frac 1{2\pi}\int
e^{-iN\theta}\Pi(e^{i\theta}x,y)\,d\theta$. The estimates for
$\Pi_N(z,w)$ are then  based on the  Boutet de
Monvel-Sj\"ostrand construction of an  oscillatory integral
parametrix for the \szego kernel:
\begin{equation}\label{oscint}\begin{array}{c}\Pi (x,y) =
S(x,y)+E(x,y)\;,\\[12pt] \mbox{with}\;\; S(x,y)=
\int_0^{\infty} e^{i t \psi(x,y)} s(x,y,t ) dt\,, \qquad
E(x,y)\in
\ccal^\infty(X \times X)\,.\end{array}
\end{equation} The amplitude has the form $s \sim \sum_{k =
0}^{\infty} t^{m -k} s_k(x,y)\in S^m(X\times X\times \R^+)$. The
phase function $\psi$ is of positive type, and as described in
\cite{BSZ2},  is given by:
\begin{equation}
\psi(z,
\theta, w, \phi) = {i} \left[1 - \frac{a(z,\bar{w})}{
\sqrt{a(z)} \sqrt{a(w)}}\; e^{i (\theta - \phi)}\right]\;,
\label{psiphase}\end{equation} where
$a\in\ccal^\infty(M\times M)$ is an almost holomorphic
extension  of the function
$a(z,\bar z):=a(z)$ on the anti-diagonal $A=\{(z,\bar
z):z\in M\}$, i.e.,\ $\dbar a$ vanishes to infinite order
along $A$. We recall from
\eqref{a} that $a(z)$ describes the Hermitian metric on $L$ in
our preferred holomorphic frame at
$z_0$, so by \eqref{preferred}, we have $
a(u)=1+|u|^2+O(|u|^3)$, and hence
\begin{equation}\label{auv} a(u,\bar v)= 1+ {u\cdot\bar v}
+O(|u|^3+|v|^3)\;.\end{equation} For further background and
notation on complex Fourier integral operators we refer to
\cite{BSZ2} and to the original paper of Boutet de Monvel and
Sj\"ostrand
\cite{BS}.

As above,  denote the $N$-th Fourier coefficient of these
operators relative to the $S^1$ action by $\Pi_N = S_N + E_N$.
Since $E$ is smooth, we have $E_N(x,y) = O(N^{-\infty})$, where
$O(N^{-\infty})$ denotes a quantity which is  uniformly
$O(N^{-k})$ on $X\times X$ for all positive $k$. Then,
$E_N(z,w)$ trivially satisfies the
remainder estimates in Theorem \ref{near-far}.

Hence it is only necessary to verify that the oscillatory integral
\begin{equation}S_N(x,y)  =   \int_0^{2\pi} e^{- i
N \theta}  S( e^{i\theta} x,y) d\theta  =  \int_0^{\infty}
\int_0^{2\pi}  e^{- i N \theta+it  \psi( e^{i\theta} x,y)}
s(e^{i\theta} x,y,t) d\theta dt \end{equation} satisfies Theorem
\ref{near-far}. This follows from an analysis of the  stationary
phase method  and remainder estimate for the rescaled parametrix
 \begin{equation}\begin{array}{l}\displaystyle S^{z_0}_N\left(
\frac{u}{\sqrt{N}}, 0;  \frac{v}{\sqrt{N}}, 0\right)
= N \int_0^{\infty} \int_0^{2\pi}
 e^{ i N \left( -\theta + t\psi\big(  \frac{u}{\sqrt{N}}, \theta;
\frac{v}{\sqrt{N}}, 0\big)\right)} s\left( \frac{u}{\sqrt{N}},
\theta; \frac{v}{\sqrt{N}}, 0, Nt\right) d\theta dt
\,,\end{array}\label{SN}\end{equation} where we changed  variables
$t \mapsto N t$. For background on the stationary phase method
when the phase is complex we refer to \cite{H}. We are
particularly interested in the dependence of the stationary phase
expansion and remainder estimate on the parameters $(u, v)$
satisfying the constraints in (i)-(ii) of Theorem \ref{near-far}.

To clarify the constraints, we recall  from  \cite{SZ2} (95) that
the \szego kernel satisfies the following far from diagonal
estimates:
\begin{equation} \label{CONSTRAINTS}  \left|\nabla^j_h \Pi_N(z,w)\right|=O(N^{-K})\qquad \mbox{for
all } \;j,  K \; \mbox{when} \;\;
d(z,w)\geq\,\frac{N^{1/6}}{\sqrt{N}} \;. \end{equation}  Hence we
may assume from now on that $z = z_0 + \frac{u}{\sqrt{N}}, w = z_0
+ \frac{v}{\sqrt{N}}$ with
\begin{equation} \label{CONSTRAINTS2} |u|+|v|\le \delta N^{1/6}
\end{equation}
for a sufficiently small constant $\delta > 0$ (see \ref{NEG}).

By \eqref{psiphase}--\eqref{auv},  the rescaled phase in \eqref{SN} has the
form:
\begin{equation}
 \wt\Psi:= t \psi\left( \frac{u}{\sqrt{N}},\theta;  \frac{
v}{\sqrt{N}}, 0\right)  -\theta  = it \left[ 1 - \frac{a\left(
\frac{u}{\sqrt{N}},  \frac{\bar v}{\sqrt{N}}\right)}{a\left(
\frac{u}{\sqrt{N}},
\frac{\bar u}{\sqrt{N}}\right)^{\half} a\left( \frac{v}{\sqrt{N}},  \frac{\bar
v}{\sqrt{N}}\right)^{\half}}\;  e^{i \theta}\right] -  \theta
\end{equation}
and the  $N$-expansion
\begin{equation}\label{entirephase} \wt\Psi=it[ 1 - e^{i \theta}]
-  \theta -\frac{it}{N}\psi_2(u,v) e^{i \theta} + t
R_3^\psi(\frac{u}{\sqrtn},\frac{v}{\sqrtn}) e^{i \theta}\,,
\end{equation}
where
\begin{equation*} \psi_2(u,v) =
u \cdot\bar{v} - \half(|u|^2 + |v|^2)=-\half|u-v|^2+i\,\Im(u\cdot\bar
v)\end{equation*} is the phase function  of (\ref{heisen-N}). After
multiplying by $iN$, we move the last two terms of \eqref{entirephase}
into the amplitude.
 Indeed, we    absorb all of  $ \exp\{(\psi_2 +i
NR_3^\psi)te^{i \theta}\}$   into the amplitude so that (\ref{SN})
is an oscillatory integral
\begin{equation}\label{phase-amplitude}  N
\int_0^{\infty} \int_0^{2\pi}
e^{iN\Psi(t,\theta)}A(t,\theta;z_0,u,v)d\theta dt + O(N^{-\infty})
\end{equation}with phase
\begin{equation}\label{phase}\Psi(t,\theta): = it ( 1 -
e^{i\theta})- \theta\end{equation} and with amplitude
\begin{equation}\label{amplitude}
A(t,\theta;z_0,u,v):=  e^{ t e^{i \theta} \psi_2(u,v)  + it e^{i
\theta} N R_3^\psi( \frac{u}{\sqrt{N}},\frac{ v}{\sqrt{N}})}\,
s\big( \frac{u}{\sqrt{N}}, \theta; \frac{v}{\sqrt{N}}, 0, Nt\big).
\end{equation}

 The  phase $\Psi$ is independent of the
parameters $(u,v)$, satisfies $\Re (i \Psi) = - t (1 - \cos
\theta) \leq 0$ and has a unique critical point at $ \{t=1, \theta
= 0\}$ where it vanishes.

The factor $e^{ t e^{i \theta} \psi_2(u,v) }$ is of exponential
growth
 in some regions.   However, since it is  a
 rescaling of a complex phase of positive type, the complex phase
 $i N \Psi +  t e^{i \theta} \psi_2(u,v)$ is of positive type,
 \begin{equation} \label{NEG} \Re (i N \Psi +  t e^{i \theta} \psi_2(u,v))
   <
 0\end{equation}
 once the cubic remainder $N t e^{i \theta}
R_3^\psi(\frac{u}{\sqrtn},\frac{v}{\sqrtn})$ is smaller than $i N
\Psi + t e^{i \theta} \psi_2(u,v)$, which occurs
 for all $(t, \theta, u, v)$ when $(u, v)$ satisfy
 (\ref{CONSTRAINTS2}) with $\delta$ sufficiently small.

 To estimate the
 joint rate of decay in $(N, u, v)$, we  follow the  stationary
phase expansion and remainder estimate in Theorem 7.7.5 of
\cite{H}, with extra attention to    the unbounded  parameter $u$.

The first step is to use a  smooth partition of unity $\{\rho_1(t,
\theta),\rho_2(t, \theta )\}$ to decompose the integral (\ref{SN})
into a region  $(1 - \epsilon, 1 + \epsilon)_t \times (-\epsilon,
\epsilon)_{\theta}$  containing the critical point and one over
the complementary set  containing no critical point. We claim that
the  $\rho_2$ integral  is of order $N^{-\infty}$ and can be
neglected. This follows by repeated partial integration as in the
standard proof together with the fact that the exponential factors
in (\ref{NEG}) decay, so that the estimates are integrable and
uniform in $u$.

We then apply  \cite{H}  Theorem 7.7.5 to  the $\rho_1$ integral.
 The first term
of the stationary phase expansion equals $N^m e^{ t e^{i \theta}
\psi_2(u,v)}$ and the remainder satisfies

\begin{equation} \label{remainder} |\wh R_{J}(P_0, u, v, N)| \leq C
N^{- m + J} \sum_{|\alpha| \leq 2J+2} \sup_{t, \theta}
|D^{\alpha}_{t, \theta} \rho_1 A(t,\theta;P_0,u,v) |.
\end{equation}
From the formula in (\ref{amplitude}) and the  fact that $s$ is a
symbol, $A$ has a polyhomogeneous expansion of the form
\begin{eqnarray} A(t,\theta;P_0,u,v) &=& \rho_1(t, \theta)   e^{ t e^{i \theta} \psi_2(u,v)}
 N^m \left[\sum_{n = 0}^{K} N^{-n/2} f_{n}(u, v;t, \theta,
P_0)+R_K(u, v, t,\theta)\right]\,,\nonumber\\
&&\qquad\qquad |\nabla^jR_{Nk}(u,v)|\le C_{jk\ep b}e^{\ep (|u|^2 +
|v|^2)} N^{-\frac{K+1}{2}}\big). \label{Alarge}\end{eqnarray} The
exponential remainder factor $e^{\ep (|u|^2 + |v|^2)}$  comes from
the fact   $\Re e^{i \theta} \psi_2 = \cos \theta \Re \psi - \sin
\theta \Im \psi$ with $\Re \psi \leq 0$ and  $|\sin \theta| <
\epsilon$ on the support of $\rho_1$. Hence, the supremum of the
amplitude in a neighborhood of the stationary phase set (in the
support of $\rho_1$) is bounded by $e^{\epsilon |\Im \psi_2|}$.
The remainder term is smaller than the main term asymptotically as
$N \to \infty$ as long as $(u,v)$ satisfies (\ref{CONSTRAINTS2}).
Part(i) of Theorem \ref{near-far} is an immediate consequence of
(\ref{Alarge}) since
$e^{\ep (|u|^2 + |v|^2)} \leq N^{\epsilon}$ for $|u| + |v|
\leq \sqrt{\log N}$.

To prove part (ii), we may assume from
(\ref{CONSTRAINTS})--(\ref{CONSTRAINTS2}) that $\sqrt{\log N} \leq
|u|+|v|  \leq\de\,N^{1/6}$. In this range the
asymptotics (\ref{Alarge}) are valid. We first rewrite the horizontal
$z$-derivatives $\frac{\d^h}{\d z_j}$ as $u_j$ derivatives, which for
$L^N$ have the form  $\sqrt{N}
\frac{\d}{\d u_j} - N A_j(\frac u{\sqrtn})$ and thus $\nabla_h$
contributes a factor of
 $\sqrtn$. We thus obtain an
asymptotic expansion and remainder for 
$\nabla^j_h
\Pi_N(z,w)$ by applying
$\nabla^j_h$ to the expansion (i) with $k = 0$: $$
\Pi^\H_1(u,\theta;v,\phi)\left[1+ N^{-1/2}
R_{N0}(u,v)\right]. $$  The operator $\nabla^j_h $ contributes a
factor of $N^{j/2}$ to each term, and thus
\begin{eqnarray*}\left|\nabla^j_h
\Pi_N(z,w)\right|& = & O\left(N^{m + j/2} \, e^{-(1 - \epsilon)
\frac{|u|^2 + |v|^2}{2}}\right)  \\
& = & O(N^{ -k })\qquad \mbox{uniformly for }\;\;
{|u|^2 + |v|^2} \ge (j+2k+2m+\ep'){ \log N} \;,
\end{eqnarray*} where $\ep'=(j+2k+2m+1)\ep$.
\qed

\section{The sharp variance estimate: Proof of Theorem
\ref{sharp}}\label{s-sharp}
We first give the proof of Theorem
\ref{sharp}, which uses the same method as  the proof of
Theorem \ref{volume}, but has simpler computations.

  We begin with some off-diagonal asymptotics for the function
$Q_N=\wt G\circ P_N$ defined in \eqref{QN}. By Proposition
\ref{DPdecay}, we see that
 \begin{equation}\label{FNfar} Q_N( z,w)\le
 \frac 1{4\pi^2} \int_0^{C/N^{2m}} \frac{-\log(1-s)}{s}\,ds
 =O\left(\frac
1{N^{2m}} \right)\;,\quad \mbox{for }\ d(z,w)>\frac{b\sqrt{\log
N}}{\sqrtn},\end{equation} with $b>\sqrt{2m}$.

Next we show the near-diagonal estimate
\begin{equation}\label{FNnear}  Q_N\Big(z_0,z_0+\frac v
 {\sqrtn}\Big)= \wt G(e^{-\half
|v|^2})+ O(N^{-1/2+\epsilon})\;,\qquad \mbox {for }\ |v|\le
b\sqrt{\log N }.\end{equation}
To verify \eqref{FNnear}, we apply Proposition \ref{best}.  Since
$P_N(z_0,z_0)=1$ and $\wt G'(t)
\to \infty$ as $t\to 1$, we need a short argument: let
\begin{equation}\label{Lambda} \Lambda_N= -\log
P_N\;.\end{equation} Recalling \eqref{Gtilde1}, we write,
\begin{equation}\label{F} F(\la):=\wt G(e^{-\la}) = -\frac
1{2\pi^2} \int_\la^\infty
\log(1-e^{-2s})\,ds \qquad\quad (\la>0)\;,\end{equation} so that
$$Q_N=F\circ \Lambda_N\;.$$
By Proposition \ref{best}, $$
\Lambda_N\left(z_0,z_0+\frac
v{\sqrtn}\right) ={\half |v|^2} + \wt R_N(v)\;,$$ where
\begin{equation}\label{RN}\wt
R_N=-\log(1+R_N)=O(|v|^2N^{-1/2+\ep})\end{equation}
By \eqref{F} and
Proposition \ref{best},
\begin{equation}\label{F'} |F'(\la)| = -\frac
1{2\pi^2}\log(1-e^{-2\la}) \le \frac
1{2\pi^2}\max\left(\log \frac 1\la,
1\right)=O\left(1+\log^+
\frac1{|v|}\right)\;.
\end{equation}
Since  $\half |v|^2+\wt R_N(v) = |v|^2\left(\half
+o(N)\right)$, it follows from
 \eqref{RN}--\eqref{F'} that
\begin{eqnarray*}  Q_N\Big(z_0,z_0+\frac v
 {\sqrtn}\Big)&=&F\left(\half |v|^2+\wt R_N(v)\right)\nonumber\\&
= & F\left(\half |v|^2\right) +O\left(\left[1+\log^+
\frac1{|v|}\right]\wt R_N(v) \right)\nonumber \\&=& \wt G(e^{-\half
|v|^2})+ O(N^{-1/2+\epsilon})\;,\qquad \mbox {for }\ |v|\le b\sqrt{\log N
},\end{eqnarray*} which gives \eqref{FNnear}.  (This computation also shows
that $Q_N$ is $\ccal^1$ and has vanishing first derivatives on the diagonal
in $M\times M$.)

Next, we note that $\wt G(t)= \frac 1{4\pi^2}\left(t^2 +\frac
{t^4}{2^2} +
\frac {t^6}{3^2}+\cdots +
\frac {t^{2n}}{n^2} +\cdots\right)$, and hence
\begin{eqnarray}\label{zeta} \int_{\C^m}\wt G(e^{-\half
|v|^2})\,dv &=& \frac 1{4\pi^2}\sum_{k=1}^\infty\int_{\C^m}
\frac{e^{-k|v|^2}}{k^2}\,dv\nonumber \\&=&  \frac
1{4\pi^2}\sum_{k=1}^\infty \frac {\pi^m}{k^{m+2}} \ = \ \frac
{\pi^{m-2}}4 \;\zeta(m+2) \;.\end{eqnarray}

 By Proposition \ref{BIPOT}, we have
\begin{equation}\label{int1s}\var\big(Z_{s_N},\phi\big)
=
\int_{M}\ical(z)\,i\ddbar\phi(z)\;,\end{equation} where
\begin{equation}\label{int2s}\ical(z)=\int_{\{z\}\times M}
Q_N(z,w)\,i\ddbar \phi(w)\;.\end{equation}

 We let
$$\Om=\frac 1{m!} \om^m$$ denote the volume form of $M$, and we write
\begin{equation}\label{psi}i\ddbar \phi =
\psi\,\Om\;,\qquad \psi \in\ccal^1.\end{equation}
To evaluate $\ical(z_0)$ at a fixed point $z_0\in M$, we choose a normal
coordinate chart centered at $z_0$ as in \S \ref{off}.
 By \eqref{FNfar} and \eqref{int2s}--\eqref{psi},
\begin{equation}\label{Iz}\ical(z_0) =   \int_{|v|\le b\sqrt{\log N}}
Q_N\left(z_0,z_0+\frac v{\sqrtn}\right)\, \psi\left(z_0+\frac
v{\sqrtn}\right)\,
\Om\left(z_0+\frac
v{\sqrtn}\right)+O\left(\frac 1 {N^{2m}}\right)\;.
\end{equation}

We recall that $\om=\frac i2 \ddbar \log a= \frac i2 \ddbar
\left[|z|^2 +O(|z|^3)\right]$ in normal coordinates.  Hence
\begin{equation}\label{OmE}\Om\left(z_0+\frac
v{\sqrtn}\right)=
\frac 1{m!}\left[ \frac i{2N} \ddbar |v|^2+ O\left(\frac
{|v|}{N^{3/2}}\right)\right]^m= \frac
1{N^m}\left[ \Om_E(v) +
 O\left(\sqrt\frac {\log N}{N}\right)\right]\;,\end{equation}
for $|v|\le b\sqrt{\log N}$, where $$\Om_E(v)= \frac 1{m!}\left(\frac
i2\ddbar |v|^2\right)^m= \prod_{j=1}^m\frac i2 dv_j\wedge d\bar v_j$$
denotes the Euclidean volume form. Since $\phi\in\ccal^3$ and hence
$\psi(z+\frac v{\sqrtn})=
\psi(z)+O(|v|/\sqrtn)$, we then have by \eqref{FNnear} and
\eqref{Iz}--\eqref{OmE},
\begin{eqnarray}\ical(z_0)&=&  \frac 1 {N^m}\left[
 \int_{|v|\le b\sqrt{\log N}}  \left\{\wt G(e^{-\half
|v|^2})+O(N^{-1/2+\ep})\right\}\,
\left\{\psi(z_0)+O(N^{-1/2+\ep})\right\}\right.\nonumber\\
&&\qquad \times \left\{ \Om_E(v) +
 O(N^{-1/2+\ep})\right\}\bigg]+O\left(\frac 1
{N^{2m}}\right)\;.\label{U1}
\end{eqnarray}
Since $\wt G(e^{-\half
|v|^2})\in \lcal^1$ by \eqref{zeta}, we have
\begin{equation}\label{U0} \ical(z_0)= \frac {\psi(z_0)} {N^m}\left[
 \int_{|v|\le b\sqrt{\log N}}  \wt G(e^{-\half
|v|^2})\Om_E(v) +
 O\left(N^{-1/2+\ep}(\log N)^m\right)\right]\;.\end{equation}

Since $\wt
G(e^{-\la})=O(e^{-2\la})$ and hence
\begin{equation}\label{U2}\int_{|v|\ge
b\sqrt{\log N}}  \wt G(e^{-\half |v|^2})\,\Om_E(v)
=O(N^{-2m})\;,\end{equation}
we can replace the integral over the
$(b\sqrt{\log N})$-ball with one over all of $\C^m$, and therefore
\begin{equation}\label{U3} \ical(z_0)= \frac {\psi(z_0)} {N^m}\left[
 \int_{C^m}  \wt G(e^{-\half
|v|^2})\Om_E(v) +
 O(N^{-1/2+\ep'})\right]= \frac {\psi(z_0)} {N^m}\left[
 \kappa_m +
 O(N^{-1/2+\ep'})\right]\;,\end{equation}
for all $\ep'>\ep$, where $\kappa_m =
 \frac{\pi^{m-2}}{4}\;\zeta(m+2)$ by \eqref{zeta}. Therefore, by
\eqref{int1s} and \eqref{U3},
\begin{equation}\label{U4} \var\big(Z_{s_N},\phi\big) = \frac 1{N^m}\int_M
\left[
 \kappa_m +
 O(N^{-1/2+\ep'})\right]\psi(z)^2\Om(z)\;.\end{equation} Since
$$\int_M \psi(z)^2\Om(z) = \int |\ddbar \phi|^2\,\Om =
\|\ddbar\phi\|^2_2\,,$$
\eqref{U4} yields the variance formula of Theorem \ref{sharp}\qed

\medskip
\section{Variance of zeros in a domain:
Proof of Theorems \ref{number}--\ref{volume}}\label{s-number}

Following the
approach of \S
\ref{s-sharp}, we now prove Theorem \ref{volume} and, as a
consequence, we also obtain Theorem
\ref{number}, which is the one-dimensional case of
Theorem \ref{volume}.

By Proposition \ref{varint2}, we have
\begin{equation}\label{int1}\var\big(\vol_{2m-2}[Z_{s_N}\cap U]\big)
=
\int_{\d U}\Upsilon
\wedge\Phi\;,\end{equation} where $\Phi=
\frac 1 {(m-1)!}\om^{m-1}$ and
\begin{equation}\label{int2}\Upsilon(z)=-\dbar_z \int_{\{z\}\times \d
U} \dbar_w Q_N(z,w)\wedge\Phi(w)\;.\end{equation}
From  \eqref{QN} and Proposition
\ref{DPdecay}, we conclude that
\begin{equation}\label{d2Qdecay}\dbar_z\dbar_wQ_N(z,w)=
O(N^{-m})\;,
 \quad \mbox{for }\ d(z,w)>b\sqrt{\frac{\log
N}N}\;,\end{equation} where we choose
$b=\sqrt{2m+3}$.  Thus, we only need to integrate \eqref{int2} over a
small ball about $z$ in $\d U$ when using
\eqref{int1}--\eqref{int2} to compute the variance.

To evaluate
$\Upsilon(z_0)$ at a fixed point $z_0\in\d U$, we choose normal
holomorphic coordinates
$\{z_j\}$ at  $z_0$, defined in a neighborhood $V$ 0f $z_0$.
By \eqref{int2}--\eqref{d2Qdecay}, we have
\begin{eqnarray}\nonumber\Upsilon(z_0)&=& - \sum_{j,k} \left(\int_{\d U}
\left.\frac{\d^2}{\d \bar z_j\d
\bar w_k}Q_N(z,w)\right|_{z=z_0}d\bar w_k\wedge
\Phi(w)\right) d\bar z_j\\ \nonumber
&=& - \sum_{j,k} \left(\int_{\left\{z_0+w\in\d
U:|w|<b\sqrt{\frac{\log N}{N}}\right\}}
\left.\frac{\d^2}{\d \bar z_j\d
\bar w_k}Q_N(z_0+z,z_0+w)\right|_{z=0}d\bar w_k\wedge
\Phi(w)\right) d\bar z_j\\&&\qquad +\ O\left(\frac
1{N^m}\right)\label{Up}\end{eqnarray}

As in the proof of Proposition \ref{sharp}, we write $Q_N=F\circ \Lambda_N$,
where $\Lambda_N= -\log
P_N$ and $F$ is given by \eqref{F}.  By Propositions
\ref{better}--\ref{best},
\begin{equation}\label{lambdaR}\Lambda_N(z_0+z,z_0+w) =\frac N2
|z-w|^2 + \wt R_N(z,w)\;,\end{equation}
where
$\wt R_N(z,w)=
-\log\left[1+R_N\left(\sqrtn\,z,\sqrtn\, w\right)\right]$ satisfies the
estimates
\begin{eqnarray}&|\wt
R_N(z,w)|=O(|z-w|^2N^{1/2+\ep}),\quad |\nabla\wt
R_N(z,w)|=O(|z-w|N^{1/2+\ep}),\nonumber \\& |\nabla^2\wt
R_N(z,w)|=O(N^{1/2+\ep}),\qquad \mbox{for }\ |z|+|w|<b\sqrt{\frac{\log
N}N}.\label{Rtilde}\end{eqnarray}  By \eqref{F},
$$ F'(\la) = \frac 1{2\pi^2}\log(1-e^{-2\la}), \qquad  F''(\la)=
\frac 1 {\pi^2(e^{2\la} -1)}\;.$$

We now let $\la=\La_N(z_0+z,z_0+w)$. By
\eqref{lambdaR}--\eqref{Rtilde},
\begin{eqnarray} |F'(\la)| & \le& \frac 1{2\pi^2}\max\left(\log \frac
1\la, 1\right)=O\left(1+\log^+
\frac1{|z-w|}\right)\;,\nonumber\\ F''(\la)&\le & \frac
1{2\pi^2\la} =O\left(\frac 1{N|z-w|^2}\right)\;,\qquad \mbox{for }\
|z|+|w|<b\sqrt{\frac{\log
N}N}\;.\label{F''}
\end{eqnarray}

Hence, for $|w|<b\sqrt{\frac{\log
N}N}$, we have
\begin{eqnarray*} &&\hspace {-.5in}\left.\frac{\d^2}{\d \bar z_j\d
\bar w_k}Q_N(z_0+z,z_0+w)\right|_{z=0}\  =\
\left.\left[F''(\la)\,\frac{\d\la}{\d\bar z_j} \,\frac{\d\la}{\d\bar
w_k} +F'(\la)\,\frac{\d^2 \la}{\d \bar z_j\d
\bar w_k}\right]\right|_{z=0}\nonumber\\&& =\
F''(\la)\left[-\half N
w_j+O(|w|N^{1/2+\ep})\right]\left[\half
Nw_k+O(|w|N^{1/2+\ep})\right]\nonumber\\&&\qquad + F'(\la)\cdot
O(N^{1/2+\ep})\nonumber
\\&& =\ -\frac 14 N^2F''(\la)
w_jw_k  +O(N^{1/2+\ep})
\left(1+\log^+\frac 1{|w|}\right)\;.\end{eqnarray*}
Furthermore, since \begin{equation}\label{F3}-F^{(3)}(t)=\half
\mbox{csch}^2 t\le t^{-2}\,,\end{equation} we have $$F''(\la)= \frac
1{\pi^2(e^{N|z-w|^2} -1)} +O(|z-w|^{-2}N^{-3/2+\ep})\;,$$ and hence
\begin{eqnarray}\nonumber &&\hspace {-.5in}\left.\frac{\d^2}{\d \bar z_j\d
\bar w_k}Q_N(z_0+z,z_0+w)\right|_{z=0}
\\&& =\ \frac{-N^2} {4\pi^2(e^{N|w|^2} -1)}
w_jw_k  +O(N^{1/2+\ep})
\left(1+\log^+\frac 1{|w|}\right)\;.
\label{double-dbar}\end{eqnarray}
We note that under our hypothesis that $\d U$ is piecewise $\ccal^2$
without cusps, we have the
estimate\begin{equation}\label{log+}\int_{\{z_0+w\in\d
U:|w|<\de\}}\left(1+\log^+\frac 1{|w|}\right)\,d\vol(w) =
O\left(\de^{2m-1}\,|\log\de|\right)\,,\quad\mbox{for }\
\de\le \half\;.\end{equation}  Substituting
\eqref{double-dbar}--\eqref{log+}  into
\eqref{Up}, we obtain
\begin{equation}\Upsilon(z_0)=
\sum_{j,k}\left(\frac {N^2}{4\pi^2} \int_{\left\{z_0+w\in\d
U:|w|<b\sqrt{\frac{\log N}{N}}\right\}}
\frac{w_jw_k d\bar w_k} {e^{N|w|^2} -1}
\wedge\Phi(w)+
O\left(N^{-m+1+2\ep}\right)\right) d\bar z_j .\label{Up1}\end{equation}

We first consider the case where $\d U$ is $\ccal^2$ smooth (without
corners). We can
 choose our normal coordinates
$\{z_j\}$ about $z_0$ so that the
real hyperplane
$\{\Im z_1=0\}$ is  tangent to $\d U$ at $z_0=0$.   We can
then write (after shrinking $V$ if necessary),
$$U\cap V=\{z\in V:\Im z_1+\phi(z)>0\}\;,$$ where  $\phi:V\to\R$ is a
$\ccal^2$ function of $(\Re\! z_1,\, z_2,\dots,z_m)$ such that
$\phi(0)=0,\ d\phi(0)=0$.

We let
$$\tau(v)=(v_1+i\phi(v),v_2,\dots, v_m)\;,$$ so that $\d U=\{\Im
v_1=0\}$ in terms of the (non-holomorphic) $v$ coordinates. We make the
change of variables
$$w=\tau_N(v):=\tau\left(\frac v{\sqrtn}\right)$$ in the integral
\eqref{Up1}:
\begin{equation}\Upsilon(z_0)=
\sum_{j,k}\left(\frac {N^2}{4\pi^2}\int_{B_N^{2m-1}}\tau_N^*\left[
\frac{w_jw_k d\bar w_k} {e^{N|w|^2} -1}
\wedge
\Phi(w)\right]+
O\left(N^{-m+1+2\ep}\right)\right) d\bar z_j ,\label{Up2}\end{equation}
where
 $$\left\{v\in\R\times\C^{m-1}:|v|<(b-1)\sqrt{\log
N}\right\} \subset B_N^{2m-1}\subset
\left\{v\in\R\times\C^{m-1}:|v|<(b+1)\sqrt{\log
N}\right\}.$$

To  evaluate the integrand in  \eqref{Up2}, we first note that
$$w_1 = \frac {v_1}{\sqrtn} + O\left(\frac{|v|^2}{N}\right)\,,\quad
d\bar w_1 = \frac 1{\sqrtn}\,d\bar v_1 + O\left(\frac{|v|}{N}\right),\
\quad w_2 = \frac {v_2}{\sqrtn},\ \dots, \ w_m = \frac
{v_m}{\sqrtn}\,.$$
Thus $N|w|^2=|v|^2+O\left(\frac {|v|^4}N\right)$, and hence by
\eqref{F3},
$$\frac 1{e^{N|w|^2}-1} = \frac 1{e^{|v|^2}-1}
+O\left(\frac{|v|^4}N\right)O(|v|^{-2}) =  \frac 1{e^{|v|^2}-1}
+O(N^{-1+\ep})\,,$$ for $ |v| <2b\sqrt{\log N}$.
Finally, we have
\begin{equation}\label{PhiE}\tau_N^*\Phi=
\frac 1{(m-1)!}\left[ \frac i{2N} \ddbar |v|^2+ O\left(\frac
{|v|}{N^{3/2}}\right)\right]^{m-1}= \frac
1{N^{m-1}}\left[ \Phi_E(v) +
 O\left(\sqrt\frac {\log N}{N}\right)\right]\end{equation}
on $B_N^{2m-1}$, where $$\Phi_E(v)= \frac 1{(m-1)!}\left(\frac i2\ddbar
|v|^2\right)^{m-1}\;.$$
Therefore, \eqref{Up2} becomes
\begin{equation}\label{IBN}\Upsilon(z_0) =  N^{-m+3/2}\sum_{j=1}^m
\left[\frac 1{4\pi^2 }\int_{B_N^{2m-1}} \frac{v_j} {e^{|v|^2} -1}
\,\dbar |v|^2 \wedge \Phi_E(v)+O(N^{-1/2\,+2\ep})\right] d\bar
z_j  \;.\end{equation}

We note that \begin{equation}\label{IBN1}
\int_{B_N^{2m-1}} \frac{v_j} {e^{|v|^2} -1}
\,\dbar |v|^2 \wedge \Phi_E(v) = \int_{B_N^{2m-1}}
\frac{v_j v_1} {e^{|v|^2} -1}\,d\vol_{\R\times
\C^{m-1}}(v)\;.\end{equation}
Since
$$\int_{x\in\R^n:|x|>b\sqrt{\log N}}\frac {|x|^2}{e^{|x|^2}-1}\,dx
=O(N^{-b^2+1})\;,$$ we can replace the integral in \eqref{IBN1} with an
affine integral, and hence
\begin{eqnarray}\label{Inu}\Upsilon(z_0) &=& N^{-m+3/2}\sum_{j=1}^m
\left[\frac 1{4\pi^2 }\int_{\R\times \C^{m-1}}\frac{v_j v_1} {e^{|v|^2} -1}\,d\vol_{\R\times
\C^{m-1}}(v)+O(N^{-1/2\,+2\ep})
\right] d\bar z_j\nonumber \\ &=&
N^{-m+3/2}\, \nu_m \;
d\bar z_1 +  O(N^{-m+1+2\ep})\;,\end{eqnarray} where
\begin{eqnarray*} \nu_m &=& \frac 1 {4                                         \pi^2} \int_{\R^{2m-1}}
\frac{x_1^2}{e^{|x|^2}-1}\,dx\ =\ \frac 1 {4\pi^2(2m-1)}
\int_{\R^{2m-1}} \frac{|x|^2}{e^{|x|^2}-1}\,dx
\\&=&  \frac 1 {4\pi^2(2m-1)} \,
\frac{2\pi^{m-1/2}}{\Gamma(m-1/2)}\int_0^\infty
\frac{r^{2m}}{e^{r^2}-1}\,dr
\\&=& \frac {\pi^{m-5/2}}{4\,\Gamma(m+1/2)}\sum_{k=1}^\infty
\int_0^\infty e^{-kr^2}\, r^{2m}\,dr
\\&=&  \frac {\pi^{m-5/2}}{4\,\Gamma(m+1/2)}\sum_{k=1}^\infty
\frac {\Gamma(m+1/2)}{2\,k^{m+1/2}} \ =\
\frac{\pi^{m-5/2}}{8}\;\zeta\Big(m+\half\Big)\;.
\end{eqnarray*}
Substituting \eqref{Inu} into \eqref{int1}, we obtain
the formula of Theorem \ref{volume}, which completes the proof for the
case where $\d U$ is smooth.

We now consider the general case where $\d U$ is piecewise smooth
(without cusps).
Let $S$ denote the set of singular points (`corners') of $\d U$, and
let $S_N$ be the small neighborhood of $S$ given by
$$S_N=\left\{z\in \d U:
\dist (z,S)<\frac{b'\sqrt{\log N}}{\sqrtn}\right\} \;,$$ where $b'>0$ is
to be chosen below. We shall show that:
\begin{itemize} \item[i) ] \eqref{Inu} holds uniformly for $z_0\in \d
U\sm S_N$;

\item[ii) ]
$ \displaystyle\sup_{z\in\d U\sm S} |\Upsilon(z)|= O\left(
N^{-m+3/2+\ep}\right)\;.$ \end{itemize}
Since $\vol_{2m-1}S_N
=O\left(\frac{\sqrt{\log N}}{\sqrtn}\right)$, the estimate
(ii) implies that the integral in  \eqref{int1} over the small set
$S_N$ is negligible and hence
\begin{equation}\label{int1a}\var\big(\vol_{2m-2}[Z_{s_N}\cap U]\big)
=\int_{\d U\sm S_N}\Upsilon
\wedge\Phi + O(N^{-m+1+2\ep})\;.\end{equation}
It then follows from (i) and \eqref{int1a} that
\begin{eqnarray*} \var\big(\vol_{2m-2}[Z_{s_N}\cap U]\big)
&=&
N^{-m+3/2}\left[\nu_m\,\vol_{2m-1}(\d U\sm S_N)
 +O(N^{-\half +2\ep})\right]\\
&=&
N^{-m+3/2}\left[\nu_m\,\vol_{2m-1}(\d U)
 +O(N^{-\half +2\ep})\right]\;,\end{eqnarray*} which is our desired
formula.

It remains to prove (i)--(ii). To verify  (ii), for each point
$z_0\in
\d U\sm S$, we choose holomorphic coordinates $\{z_j\}$ and
non-holomorphic coordinates $\{v_j\}$ as above.  We can choose these
coordinates on a geodesic ball $V_{z_0}$ about $z_0$ of a fixed radius
$R>0$ independent of the point $z_0$, but if $z_0$ is near a corner,
$\d U$ will coincide with $\{\Im v_1=0\}$ only in a small neighborhood
of
$z_0$. To be precise, we let $K_{z_0}$ denote the connected component of
$V_{z_0}\cap \d U\sm S$ containing $z_0$. Then we choose
$\phi \in\ccal^2(V_{z_0})$ with
$\phi(0)=0,\ d\phi(0)=0$, such that
\begin{equation}\label{Kz}\{z\in V:\Im z_1+\phi(z)=0\}\supset
K_{z_0}\;.\end{equation}

Choose $N_0>0$ such that  $b'\sqrt{\frac{\log N_0}{N_0}} <R$; then
$$\textstyle\left\{w\in\d
U:d(z_0,w)<b'\sqrt{\frac{\log N}{N}}\right\}\subset V_{z_0}\;,\quad
\mbox{for }\ N\ge N_0\;.$$
Then for $N\ge N_0$,  the integrals
\eqref{Up2} and
\eqref{IBN} hold, except now they are over a piecewise smooth
hypersurface
$\wt B_N$ of the
$(b'\sqrt{\log N})$-ball in
$\C^m$ instead of the linear hypersurface $B_N^{2m-1}$.
Since $$\left| \frac{v_j\,\dbar
|v|^2} {e^{|v|^2} -1}  \right| \le  \frac{|v|^2} {e^{|v|^2} -1}\le 1\;,$$ where the above
norm is respect to the Euclidean metric on $T^*(\C^m_{\{v\}})$,
it follows from \eqref{IBN} (with $B_N$ replaced by $\wt B_N$) that
\begin{equation}\label{upper}|\Upsilon(z_0)|
\le \frac {m+o(1)}{4\pi^2}\, N^{-m+3/2}\,\vol_{2m-1}^E(\wt B_N)
\;,\end{equation} where $\vol^E$ denotes Euclidean volume. Since
$\d U$ is piecewise smooth, we see that
\begin{eqnarray}\label{volE}\nonumber\vol_{2m-1}^E(\wt
B_N)&\le&\textstyle N^{m-1/2}\,[1+o(1)]\,\vol_{2m-1}\left\{w\in\d
U:d(z_0,w)<b'\sqrt{\frac{\log N}{N}}\right\}\\& = & \textstyle
N^{m-1/2}\; O\left(\left[\frac{\log
N}{N}\right]^{m-1/2}\right)\   =\ O\left((\log
N)^{m-1/2}\right)\;.\end{eqnarray} Combining
\eqref{upper}--\eqref{volE}, we obtain the bound (ii).

To verify (i), we let $$C= \sup _{z\in \d U\sm S}\;
\frac{\dist (z,S)}{\dist (z,\d U\sm K_z)}\;.$$  We recall that our
assumption that $\d U$ `has no cusps' means that $\overline U$ is
locally $\ccal^2$ diffeomorphic to a polyhedral cone, which implies
that  $C<+\infty$.  We now let $b'=Cb$, where $b=\sqrt{2m+3}$ as before.

Consider any point $z_0\in \d U\sm S_N$, $N\ge N_0$.  Then
$$\dist (z_0,\d U\sm K_{z_0})\ge \frac {\dist (z,S)}C \ge
\frac{b'\sqrt{\log N}}{C\,\sqrtn} = \frac{b\sqrt{\log N}}{\sqrtn}\;.$$
Thus by our far-off-diagonal decay estimate \eqref{d2Qdecay}, the
points in $\d U\sm K_{z_0}$ contribute negligibly to the integral in
\eqref{Up2}, so that integral can be taken over $\tau_N\inv(K_{z_0})$,
or over the linear $(b\sqrt{\log N})$-ball $B_N$. Then
\eqref{Inu} follows as before.

Thus we have verified (i)--(ii), which
completes the proof of Theorem \ref{volume} for the general case
where $\d U$ has corners.\qed

\section{Asymptotic normality: Proof of Theorem
\ref{AN}}\label{s-normality}

The proof is a combination of  Propositions \ref{DPdecay} and
\ref{better} with a general result of Sodin-Tsirelson \cite{ST} on
asymptotic normality of nonlinear functionals of Gaussian
processes. Following \cite{ST}, we define a {\it normalized complex
Gaussian process\/} to be a random function $w(t)$ on a measure
space $(T,
\mu)$ of the form
$$w(t)=\sum c_j  f_j(t)\;,$$ where the $c_j$ are i.i.d.\ complex
Gaussian random variables (of mean 0, variance 1), and the $f_j$ are
(fixed) complex-valued measurable functions such that
$$\sum|f_j(t)|^2=1\quad
\mbox{for all }\ t\in T.$$ 
We let $w_1,w_2,w_3,\dots$ be a sequence of normalized complex
Gaussian processes on a finite measure space $(T, \mu)$. Let $f(r)
\in \lcal^2(\R^+, e^{-r^2/2} rdr)$ and let
$\psi:  T
\to
\R$ be bounded measurable. We write
$$Z_N^{\psi}(w_N) = \int_T f(|w_N(t)|) \psi(t) d\mu(t).$$

\begin{theo} \cite[Theorem 2.2]{ST} Let $\rho_N(s, t)$ be the
covariance functions for the Gaussian processes $w_N(t)$. Suppose that
\begin{enumerate}

\item[i)]  \quad $\displaystyle\liminf_{N \to \infty} \frac{\int_T
\int_{T} |\rho_N(s, t)|^{2 \alpha} \psi(s) \psi(t) d\mu(s)
d\mu(t)}{\sup_{ s \in T} \int_T |\rho_N(s, t)| d\mu(t)} > 0\;,$\\ for
$\al=1$ if  $f$ is monotonically increasing, or for all $\al\in\Z^+$
otherwise;

\item[ii)]  \quad $\displaystyle\lim_{N \to \infty}\; \sup_{s \in T}
\int_T |\rho_N(s, t)| d\mu(t) = 0.$

\end{enumerate}
Then the distributions of the random variables
$$ \frac{Z_N^{\psi} - \E Z_N^{\psi}}{\sqrt{\var(Z_N^{\psi})}}$$
converge weakly to
$\ncal(0, 1)$ as $ N \to \infty$.
\end{theo}

We apply this result in the case $f(r)=\log r$ and $(T, \mu) = (M,
\Om)$, with the normalized Gaussian processes 
$$w_N(z) : = \frac{s_N(z)}{\sqrt{\Pi_N(z,z)}}, $$
where $s_N$ is a random
holomorphic section in $H^0(M, L^N)$ with respect to its Hermitian
Gaussian measure.  The covariance kernel of this Gaussian process is 
$P_N(z,w)$. Further, we let $\phi$ be a  $\ccal^3$ real  $(2m-2)$-form and we
write $\ddbar
\phi =
\psi
\Omega$ as before (and hence $\psi \in \ccal^1$), so that
$$Z_N^{\psi}(w_N) =(Z_{s_N}, \phi) = \int_{M} \log |s_N(z)| \ddbar
\phi(z)$$ is the smooth linear statistic of
integration of the fixed test form $\phi$  over the random zero set.
(This was the
application of interest in \cite{ST}, where they considered random
functions on $\C$, $\CP^1$, and the disk.)

By Proposition \ref{EZ}, we have
$$\E Z_N^{\psi}  = \frac{\sqrt{-1}}{2\pi}\ddbar \log\Pi_N(z,z)+N\om, $$
hence
$$Z_N^{\psi}(w_N) - \E Z_N^{\psi} =
\int_{M} \log \frac{|s_N(z)|}{\sqrt{\Pi_N(z,z)}} \ddbar
\phi(z),\;\; \phi \in \dcal^{m-1, m-1}(M). $$

To apply the theorem it suffices  to check that $P_N(z,w)$
satisfies conditions (i)--(ii). We start with (ii): by Proposition
\ref{DPdecay},
$$\lim_{N \to \infty}\; \sup_{z \in M} \int_{d(z,w) > \sqrt{\frac
{b\log N}{N}}} |P_N(z,w)| dV_{\omega}(w)= 0. $$ On the other hand,
since $|P_N(z,w)| \leq 1$ it is obvious that the same limit holds
for $d(z,w) \leq \sqrt{\frac{b \log N}{N}}.$

 To check (i), we again break up the integral into the near
 diagonal $d(z,w) \leq \sqrt{\frac{b \log N}{N}}$ and the
 off-diagonal $d(z,w) > \sqrt{\frac{b \log N}{N}}$. As before, the
 integrals over the off-diagonal set tend to zero rapidly and can
 be ignored in both the numerator and denominator.

On the near diagonal, we  replace $P_N$ by its asymptotics in
Proposition \ref{better}. The asymptotics are constant in $z$ and
with uniform remainders, so the condition becomes
$$\liminf_{N \to \infty} \frac{\int_M \int_{|u| < \sqrt{b\log N} }
e^{- |u|^2}|1 + R_N(u)|^{2 } \psi(z +
\frac{u}{\sqrt{N}}) \psi(z)\, du\, dV_{\omega}(z)}{ \int_{|u| <
\sqrt{b\log N} } e^{-\half |u|^2}|1 + R_N(u)|\, du} >
0. $$ Since $\psi \in \ccal^1$, the ratio clearly tends to $2^{-m}\int_M
\psi(z)^2 dV_{\omega} > 0$, completing the proof.\qed


\begin{thebibliography}{WWW}

\bibitem[Ber]{Ber} F. A. Berezin,
Quantization,  {\it Izv.\ Akad.\ Nauk SSSR Ser.\ Mat.} 38 (1974),
1116--1175 (Russian). English translation: {\it Math.\ USSR-Izv.} 38
(1974), no. 5, 1109--1165 (1975).

\bibitem[Be]{Be} B. Berndtsson, Bergman kernels related to Hermitian line
bundles over compact complex manifolds,  {\it Explorations in Complex and
Riemannian Geometry: A Volume Dedicated to Robert E. Greene\/}, Contemporary
Math., vol.\ 332, Amer.\ Math.\ Soc., Providence, RI, 2003.


\bibitem[BD]{BD} P. Bleher and X. Di,  Correlations between zeros of a
random
polynomial, {\it J. Statist.\ Phys.} 88 (1997), 269--305.


\bibitem[BR]{BR} P. Bleher and R.   Ridzal, ${\rm SU}(1,1)$ random
polynomials, {\it
J. Statist.\ Phys.} 106 (2002),  147--171.

\bibitem[BSZ1]{BSZ1} P. Bleher, B. Shiffman and S. Zelditch,
Poincar\'e-Lelong approach to universality and scaling of correlations
between zeros, {\it Comm.\ Math.\ Phys.} 208 (2000), 771--785.


\bibitem[BSZ2]{BSZ2} P. Bleher, B. Shiffman and S. Zelditch, Universality
and scaling of correlations between zeros on complex manifolds,  {\it
Invent.\
Math.} 142 (2000), 351--395.

\bibitem[BSZ3]{BSZ0} P. Bleher, B.  Shiffman, and S. Zelditch,
Universality and scaling of zeros on symplectic manifolds. {\it
Random matrix models and their applications}, 31--69, Math. Sci.
Res. Inst. Publ., 40, Cambridge Univ. Press, Cambridge, 2001.

\bibitem[BBL]{BBL} E.  Bogomolny, O.  Bohigas, and P.  Leboeuf, Quantum
chaotic dynamics and random polynomials, {\it J.  Statist.\ Phys.} 85 (1996),
639--679.

\bibitem[BoSj]{BS} L. Boutet de Monvel and J. Sj\"ostrand, Sur la
singularit\'e des noyaux de Bergman et de Szeg\"o, {\it
Asterisque\/} 34--35 (1976), 123--164.


\bibitem[Ch]{Ch} M. Christ, On the $\bar{\partial}$ equation in weighted $L^2$
norms in $\C^1$, {\it Jour.\ Geom.\ Anal.} 3 (1991), 193--230.

\bibitem[DLM]{DLM}
X. Dai, K. Liu and X. Ma, On the asymptotic expansion of Bergman
kernel (math.DG/0404494).


\bibitem[De]{D} H. Delin, Pointwise estimates for the weighted Bergman
projection kernel in $\mathbf C\sp n$,
using a weighted $L\sp 2$ estimate for the $\overline\partial$
equation, {\it Ann.\ Inst.\ Fourier (Grenoble)} 48 (1998),
967--997.

\bibitem[FH]{FH} P. J. Forrester and G. Honner, Exact statistical
properties
of the zeros of complex random polynomials, {\it J. Phys.\ A} 32  (1999),
2961--2981.

\bibitem[GJ]{GJ}  J. Glimm and A. Jaffe, {\it Quantum physics. A functional integral point of view.}
 Second edition. Springer-Verlag, New York, 1987.

\bibitem[Han]{Han} J. H. Hannay, Chaotic analytic zero points: exact
statistics for those
of a random spin state, {\it J. Phys.\ A} 29 (1996), L101--L105.

\bibitem[H\"o]{H} L. H\"ormander, {\it The Analysis of Linear Partial
Differential Operators, I},  Springer Verlag, N.Y.,
1983].

\bibitem[Ja]{J} S. Janson, {\it Gaussian Hilbert Spaces}. Cambridge
Tracts in Mathematics, 129. Cambridge University Press, Cambridge, 1997.

\bibitem[Kac]{Kac} M. Kac, On the average number of real roots of a random
algebraic
equation, II, {\it Proc.\ London Math. Soc.} 50 (1949), 390--408.

\bibitem[Kah]{Kah} J.-P. Kahane,  {\it Some random series of
functions\/}, Second edition, Cambridge Studies in Advanced
Mathematics, 5, Cambridge University Press, Cambridge, 1985.

\bibitem[Kn]{K} A. W. Knapp,{\it Representation theory of semisimple
groups. An overview based on examples.}
 Reprint of the 1986 original. Princeton Landmarks in Mathematics. Princeton University Press, Princeton, NJ,
 2001.

\bibitem[Li]{Li} N. Lindholm,  Sampling in weighted $L^p$ spaces of entire
functions in $\C^n$ and estimates of the Bergman kernel, {\it J. Funct.\
Anal.} 182 (2001), 390--426.


\bibitem[MM]{MaMa}  X.  Ma and  G. Marinescu,
Generalized Bergman kernels on symplectic manifolds
(math.DG/0411559).


\bibitem[NV]{NV} S.  Nonnenmacher and A.  Voros, Chaotic eigenfunctions in
phase space,  {\it J.
Statist.\ Phys.} 92 (1998), 431--518.



\bibitem[SZ1]{SZ} B. Shiffman and S. Zelditch, Distribution of zeros of
random and quantum chaotic sections of positive line bundles,
{\it Comm.\ Math.\ Phys.} 200 (1999), 661--683.

\bibitem[SZ2]{SZ2} B. Shiffman and S. Zelditch,  Asymptotics of almost
holomorphic sections of ample line bundles on symplectic
manifolds,  {\it J. Reine Angew.\ Math.}  544 (2002), 181--222.

\bibitem[SZ3]{SZ3} B. Shiffman and S. Zelditch,
Complex zeros of  random fewnomials, in preparation.

\bibitem[So]{So} M. Sodin,  Zeros of Gaussian analytic functions, {\it
Math.\ Res.\
Lett.} 7 (2000), 371--381.


\bibitem[ST]{ST} M. Sodin and B. Tsirelson,
Random complex zeros, I. Asymptotic normality, {\it Israel J. Math.} 144
(2004), 125--149.




\bibitem[Ze]{Z} S. Zelditch, \szego kernels and a theorem of Tian,
{\it Internat.\ Math.\ Res.\ Notices} 1998 (1998),  317--331.

\end{thebibliography}
\end{document}